\documentclass[11pt,leqno,tbtags]{article}
\usepackage{amsthm,amsmath,amssymb,xspace,epsfig,latexsym}
\numberwithin{equation}{section}

\newtheorem{theorem}{Theorem}[section]
\newtheorem{corollary}[theorem]{Corollary}
\newtheorem{lemma}[theorem]{Lemma}
\newtheorem{proposition}[theorem]{Proposition}
\newtheorem{claim}[theorem]{Claim}

\theoremstyle{definition}
\newtheorem{remark}[theorem]{Remark}
\newtheorem{example}[theorem]{Example}
\newtheorem{algorithm}[theorem]{Algorithm}

\newcommand{\stirlingone}[2]{%
  \begin{bmatrix}
    #1\\#2
  \end{bmatrix}%
}

\newcommand{\EE}{{\bf  E}}

\newcommand{\PP}{{\bf  P}}

\newcommand{\mm}{{\bf  m}}
\newcommand{\nn}{{\bf  n}}

\newcommand{\sss}{{\bf s}}
\newcommand{\ttt}{{\bf t}}

\newcommand{\ct}{{\tilde{c}}}
\newcommand{\dt}{{\tilde{d}}}
\newcommand{\xt}{{\tilde{x}}}

\newcommand{\Pc}{{\cal P}}

\newcommand{\hh}{\hat{h}}

\newcommand{\Kh}{\widehat{K}}

\newcommand{\xh}{\hat{x}}

\newcommand{\kappah}{\hat{\kappa}}

\newcommand{\begp}{\begin{proposition}}
\newcommand{\enp}{\end{proposition}}
\newcommand{\begt}{\begin{theorem}}
\newcommand{\ent}{\end{theorem}}
\newcommand{\begl}{\begin{lemma}}
\newcommand{\enl}{\end{lemma}}
\newcommand{\begc}{\begin{corollary}}
\newcommand{\enc}{\end{corollary}}
\newcommand{\begcl}{\begin{claim}}
\newcommand{\encl}{\end{claim}}
\newcommand{\begr}{\begin{remark}}
\newcommand{\enr}{\end{remark}}
\newcommand{\begal}{\begin{algorithm}}
\newcommand{\enal}{\end{algorithm}}
\newcommand{\begd}{\begin{definition}}
\newcommand{\enf}{\end{definition}}
\newcommand{\begx}{\begin{example}}
\newcommand{\enx}{\end{example}}
\newcommand{\bega}{\begin{array}}
\newcommand{\ena}{\end{array}}

\newcommand{\sfrac}[2]{{\textstyle\frac{#1}{#2}}}

\def\rompar(#1){\textup(#1\textup)}    

\newcommand\eps{\varepsilon}

\newcommand\gl{\lambda}

\newcommand{\refS}[1]{Section~\ref{#1}}
\newcommand{\refT}[1]{Theorem~\ref{#1}}
\newcommand{\refC}[1]{Corollary~\ref{#1}}
\newcommand{\refL}[1]{Lemma~\ref{#1}}
\newcommand{\refP}[1]{Proposition~\ref{#1}}
\newcommand{\refE}[1]{Example~\ref{#1}}
\newcommand{\refR}[1]{Remark~\ref{#1}}

\newcommand\ie{i.e.\spacefactor=1000}

\renewcommand\Re{\operatorname{Re}}

\newcommand\nopf{\qed}   
\newcommand\noqed{\renewcommand{\qed}{}} 

\newcommand{\bege}{\begin{equation}}
\renewcommand{\Pr}{\PP}

\newcommand{\ignore}[1]{}

\begin{document}

\setcounter{page}{0}
\thispagestyle{empty}

\begin{center}
{\Large \bf
Asymptotic Analysis via Mellin Transforms for Small Deviations in
$L^2$-norm of Integrated Brownian Sheets\\
}
\normalsize

\vspace{4ex}
{\sc James Allen Fill\footnotemark} \\
\vspace{.1in}
Department of Mathematical Sciences \\
\vspace{.1in}
The Johns Hopkins University \\
\vspace{.1in}
{\tt jimfill@jhu.edu} and {\tt http://www.mts.jhu.edu/\~{}fill/} \\
\vspace{.2in}
{\sc and} \\
\vspace{.1in}
{\sc Fred Torcaso}\\ 
\vspace{.1in}
Department of Mathematical Sciences \\
\vspace{.1in}
The Johns Hopkins University \\
\vspace{.1in}
{\tt torcaso@mts.jhu.edu} and {\tt http://www.mts.jhu.edu/\~{}torcaso/} \\
\end{center}
\vspace{3ex}

\begin{center}
{\sl ABSTRACT} \\
\end{center}

We use Mellin transforms to compute a full asymptotic expansion for the tail
of the Laplace transform of the squared $L^2$-norm of any multiply-integrated
Brownian sheet.  Through reversion we obtain corresponding
strong small-deviation estimates.
\bigskip
\bigskip

\begin{small}
\par\noindent
{\em AMS\/} 2000 {\em subject classifications.\/}  Primary 60G15, 41A60;
secondary 60E10, 44A15, 41A27.
%

\medskip
\par\noindent
\emph{Date.} September~23, 2003.
\end{small}

\footnotetext[1]{Research supported by NSF grant DMS--0104167,
and by The Johns Hopkins University's Acheson J.~Duncan Fund for the
Advancement of Research in Statistics.}

\newpage
\addtolength{\topmargin}{+0.5in}

\newcommand{\tab}[0]{\hspace{.1in}}

\newtheorem*{acknowledgment}{Acknowledgment}

\newpage


\section{Introduction}
\label{S:intro}

Throughout, for readability, we often use shorthand vector notation such
as~$\ttt$ for $(t_1, \dots, t_d)$ and $d\ttt$ for $dt_1 \cdots dt_d$.

A~$d$-parameter (standard) \emph{Brownian sheet} $B=\left(
B(\ttt)\!:\!\ttt \in [0, 1]^d \right)$ is defined to be a
real-valued Gaussian random field with continuous sample paths,
mean function $\EE\,B(\ttt) = 0$, and covariance kernel~$K_{\bf
0}(\sss, \ttt) := \EE[B(\sss) B(\ttt)] = \prod_{j = 1}^d
\min\{s_j,t_j\}$ for~$\sss, \ttt\in [0,1]^d$.  Notice that the
covariance operator with kernel~$K_{\bf 0}$ is the tensor product
of~$d$ copies of the covariance operator of Brownian motion.  Let
$\mm = (m_1, \dots, m_d)$ be a fixed~$d$-vector with nonnegative
integer components.  We define \emph{$\mm$-integrated Brownian
sheet}~$X_{\mm}$ by
\begin{equation}
\label{ibs}
X_{\mm}(\ttt) := \int_0^{t_1}\!\!\!\! \cdots \int_0^{t_d}\!\prod_{j = 1}^d
\frac{(t_j - u_j)^{m_j}}{m_j!}\,B(du_1, \dots, du_d).
\end{equation}
It follows immediately that~$X_{\mm}$ is a mean-zero Gaussian random field
on~$[0, 1]^d$ with covariance kernel
$$
K_{\mm}(\sss, \ttt)=\prod_{j = 1}^d \int_0^{\min\{s_j, t_j\}}\!
\frac{(s_j - u_j)^{m_j} (t_j - u_j)^{m_j}}{(m_j!)^2}\,du_j,
$$
and so the covariance operator of $X_{\mm}$ is simply the
tensor product of the covariance operators with kernels~$K_{m_j}$,
which are the covariance kernels of $m_j$-integrated Brownian
motions.  Observe that $X_{\bf 0}$ is simply $d$-parameter
Brownian sheet~$B$.

\begin{remark}
\label{R:mintBM}
To motivate the definition of an integrated Brownian sheet, recall that
$m$-times integrated Brownian motion~$X_m$ can be defined (naturally enough) in
terms of Brownian motion~$B$ by
\begin{equation}
\label{IBMdef}
X_m(t) := \int_0^t\!\int_0^{s_m}\!\!\!\!\cdots\int_0^{s_2}\!B(s_1)\,ds_1 \cdots
ds_m
\end{equation}
for integer $m \ge 1$, with $X_0 := B$.
It is not difficult to see that $X_m$ has the same distribution as the process
with value at time~$t$ given by
\begin{equation}
\label{IBMaltdef}
\frac{1}{m!} \int_0^t\!(t-u)^m\,dB(u).
\end{equation}
That is, the~$m$ integrations in~\eqref{IBMdef} can be collapsed to the one
in~\eqref{IBMaltdef}.  Similar motivation can be given for our definition of
$\mm$-integrated Brownian sheet.
\end{remark}

Consider the squared~$L^2$-norm of~$X_{\mm}$:
\begin{equation}
\label{Vdef}
V^2 \equiv V_{\mm}^2 := \int_0^1\!\! \cdots \int_0^1\!X^2_{\mm}(\ttt)\,d\ttt.
\end{equation}
The classical Karhunen--Lo\`eve expansion tells us that~$X_{\mm}$ has the same
distribution as the process $\left( \sum_{\nn} \sqrt{a_{\nn}}
\varphi_{\nn}(\ttt) \xi_{\nn} \right)$, where the $\xi_{\nn}$'s are i.i.d.\
standard normal random variables and the $\varphi_{\nn}$'s form a complete
orthonormal system of eigenvectors, with corresponding eigenvalues~$a_{\nn}$,
for the covariance operator ${\cal A}_{\mm}:L^2([0,1]^d) \to L^2([0,1]^d)$.
The spectrum $\sigma({\cal A}_{\mm}) =
\{ a_{\nn} \}$
is the product of spectra of the covariance
operators ${\cal A}_{m_j}$ of the associated $m_j$-integrated Brownian motions:
$$
\sigma({\cal A}_{\mm}) = \sigma({\cal A}_{m_1}) \otimes \cdots
\otimes \sigma({\cal A}_{m_d}),
$$
where $\otimes$ represents elementwise set multiplication:\ $S
\otimes T := \{s t:s \in S\mbox{\ and\ }t \in T\}$. It therefore follows that,
in distribution,
\begin{equation}
\label{Vsum}
V^2 = \sum_{n_1 = 1}^{\infty} \cdots \sum_{n_d = 1}^{\infty} a_{\nn}
\xi^2_{\nn},
\end{equation}
where~$a_{\nn}= a_{\nn}(\mm) = a_{n_1}(m_1) \cdots a_{n_d}(m_d)$
and~$(a_{n_j}(m_j))$ are the eigenvalues for ${\cal A}_{m_j}$.

We are interested in deriving a ``strong'' small-deviations result, \ie, in
computing the lead-order asymptotics as $\eps \to 0+$ for the small-ball
probability $\Pr(V^2 \le \eps)$ (not just for its logarithm).  In light of the
representation~\eqref{Vsum}, this can be done using Sytaja's Tauberian
Theorem~\cite{Sytaja1974}:
\begin{theorem}[Sytaja~\cite{Sytaja1974}]
\label{T:sytaja}
Suppose that $a_n > 0$ for all $n$ with $\sum_{n=1}^{\infty} a_n < \infty$, and
that $(\xi_n)_{n \geq 1}$ are i.i.d.\ standard normal random variables.  Then
as $\varepsilon \to 0+$,
$$
\Pr\left(\sum_{n=1}^{\infty} a_n \xi_n^2 \le \varepsilon \right) =
(1 + o(1)) \left[ -2 \pi (x^*)^2 h''(x^*) \right]^{-1/2} \exp\left\{- [h(x^*) -
\eps\,x^*]\right\}
$$
where $h(x) := -\log \EE\,\exp\left\{- x \sum_{n = 1}^{\infty} a_n
\xi_n^2\right\}$, $x \ge 0$, denotes the log Laplace transform and~$x^* \equiv
x^*(\eps)$ is defined implicity in terms of~$\varepsilon$:
\begin{equation}\label{xepsrel}
h'(x^*) = \varepsilon.~\qed
\end{equation}
\end{theorem}

In this paper we first show how to derive, for arbitrary~$\mm$, a complete
asymptotic expansion as $x\to \infty$ for the log Laplace transform of~$V^2$,
namely, $h_{\mm}(x) := - \log \EE\,\exp\{-x V_{\mm}^2)$, and for each of its
derivatives.  Our method proceeds in two steps.  First we study carefully the
one-dimensional functions $h_{m_j}$.
Then
we apply Mellin-transform techniques that prove to be quite powerful
in handling our ``tensored'' processes.  Although Mellin transforms enjoy
common use in the analysis of algorithms~\cite{FS2003}, their application here
is, to our knowledge, the first in studying small or large deviations.
As the reader will see, their use is very natural and makes many computations
completely transparent.

We also show how one can develop the explicit asymptotic behavior
of~$x^*$ in terms of~$\eps$ via the relation~\eqref{xepsrel}.
This ``reversion'' can be quite nontrivial and is the subject of
\refS{S:reversion}.  Once this reversion is understood, one can
then easily apply \refT{T:sytaja} to obtain the strong small-ball asymptotics.

If~$X_m$ is $m$-times integrated Brownian motion, then the lead-order
asymptotics of~$ - \log \Pr(V^2 \le \eps)$ have been studied by Chen and
Li~\cite{ChenLi2003}.  In this one-dimensional case, lead-order asymptotics for
the small-ball probability $\Pr(V^2 \le \varepsilon)$ itself have been studied
by Gao et al.~(\cite{GaoHannigTorcaso2003},~\cite{GaoHannigLeeTorcaso2003a},
and~\cite{GaoHannigLeeTorcaso2003b}), Nazarov~\cite{Nazarov2003}, and
Nazarov and Nikitin~\cite{NazarovNikitin2003}.
A closed-form expression for the
Laplace transform of $V_m^2$ for any~$m$ has been derived in Gao et
al.~\cite{GaoHannigLeeTorcaso2003b}; we should also mention that in the
special case $m = 1$ the Laplace transform of~\eqref{Vdef} had been
obtained explicitly in Khoshnevisan and Shi~\cite{KhoshnevisanShi1998} several
years earlier.

We now turn our attention to higher dimensions.  For $d \ge 1$,
Csaki~\cite{Csaki1982} computes the lead-order asymptotics of~$-\log \Pr(V_{\bf
0}^2\le \eps)$, that is, the lead order for logarithmic small deviations of
$d$-dimensional Brownian sheet.  But a different approach is needed for nonzero
$m_j$'s. To illustrate, consider~$(m_1, m_2)$-integrated Brownian sheet.
Then
\begin{eqnarray}
h(x)
    & = & - \log \EE\,\exp\left\{- x \sum_{n_1 = 1}^{\infty} \sum_{n_2 =
            1}^{\infty} a_{n_1, n_2} \xi_{n_1, n_2}^2\right\} \nonumber\\
\label{hexpression}
    & = & \frac12 \sum_{n_1 = 1}^{\infty} \sum_{n_2 = 1}^{\infty} \log
\left(1 + \frac {2x}{\kappa_{n_1}(m_1) \kappa_{n_2}(m_2)}\right),
\end{eqnarray}
where $a_{n_1, n_2} = 1/\left[\kappa_{n_1}(m_1)\kappa_{n_2}(m_2)\right]$
and the $\kappa_n(m)$, $n \geq 1$, are
the reciprocals of the eigenvalues of~${\cal A}_m$.  The reason
for moving to reciprocals is that the~$\kappa$'s
have nice representations as the zeros of certain entire
functions~\cite{GaoHannigLeeTorcaso2003b}.
If $m_1 = m_2 = 0$, then $\kappa_n(0) = \left[ \left( n - \frac{1}{2} \right)
\pi \right]^2$
and Csaki's method works well because there are known
expressions and bounds for the summations involved.  However, if a component
of~$\mm$ differs from~$0$ then his approach is untenable.  For example, if $m_1
= 0$ and $m_2 = 1$, then $(\kappa_n(1))$ are the solutions to the equation
$\cos(z^{1/4}) \cosh(z^{1/4}) + 1 = 0$.  When~$n$ is large, $\kappa_n(1)$ is
quite close to $\left[ \left( n - \frac{1}{2} \right) \pi \right]^4$ (see,
for instance,~\cite{GaoHannigTorcaso2003}), but it does not seem possible to
obtain explicit expressions. In this case we immediately see the complexity
introduced when $m_2 = 1$.

Csaki handled only the case $\mm = {\bf 0}$.  Recently, Karol' et
al.~\cite{KarolNazarovNikitin2003} extended the classical result
of Csaki~\cite{Csaki1982} to obtain lead-order logarithmic small
deviations for many Gaussian random fields; in particular, their
treatment can handle our problem for arbitrary~$\mm$.  However,
their methods do not seem to extend to obtaining asymptotic
results beyond the lead order for the log small-ball
probability.

We note that the methods developed in this paper can be applied to
many other Gaussian random fields besides $\mm$-integrated Brownian
sheet~$X_{\mm}$.  Consider a Gaussian random field whose covariance
operator is the tensor product of marginal operators.  If the Gaussian processes
corresponding to these marginal operators each have a log Laplace transform
whose asymptotic behavior can be established, then the methods of this paper
will apply.  Further, since (at least for~$X_{\mm}$) we can derive a complete
asymptotic expansion for the log Laplace transform, we need only produce a
complete asymptotic expansion extending Sytaja's theorem in order to obtain a
complete asymptotic expansion---not just the lead-order asymptotics---for the
small-ball probability.  We will produce the needed generalization of Sytaja's
theorem, and apply it to small-ball probabilities for various Gaussian random
fields, in future work.

This paper is organized as follows.  In \refS{S:prelims} we provide basic
background concerning the Mellin transform.  We derive asymptotic expansions for
$h(x)$ and its derivatives in \refS{S:hasy}.  In \refS{S:reversion} we discuss
explicit expansions for the solution $x^* \equiv x^*(\eps)$ to the
equation $h'(x) = \eps$ [recall~\eqref{xepsrel}].  Finally, in \refS{S:small} we
derive lead-order asymptotics for the small-ball probability.

Without loss of generality we will assume throughout this paper that
$$
0 \leq m_1\le m_2 \le \cdots \le m_d < \infty.
$$
For multibranched functions of a complex variable, we use the principal
branch unless otherwise specified.

\section{Mellin transforms}
\label{S:prelims}

In this section we collect, for the reader's convenience, some useful
facts concerning Mellin
transforms.  Our treatment follows closely the superb book manuscript
by Flajolet
and Sedgewick~\cite{FS2003}.  Throughout, we use the following
notation for an open vertical
strip in the complex plane:\ for real $c, d$,
$$
\langle c, d \rangle := \{s \in {\mathbb C}:c < \Re s < d\}.
$$

\subsection{Definition}
\label{S:mtdef}

Suppose that $f:[0,\infty)\to {\mathbb C}$ is locally integrable
[\ie,\ is Lebesgue integrable over any bounded closed subinterval of
$(0, \infty)$;
this condition is met, for example, if~$f$ is continuous].
For those $s \in {\mathbb C}$ such that $x \mapsto f(x) x^{s - 1}$ is
Lebesgue integrable
over $(0, \infty)$, we define the \emph{Mellin transform} $f^*$ of
$f$ as the Lebesgue
integral
\begin{equation}\label{mellin}
f^*(s) := \int_0^{\infty}\!f(x) x^{s - 1}\,dx.
\end{equation}
Keep in mind that Lebesgue integrability is a form of absolute integrability.

\subsection{Existence:\ the fundamental strip; transform of derivative}
\label{S:mtexist}

Given local integrability, what is at issue for the existence
of~\eqref{mellin} is the
behavior of~$f$ near~$0$ and near~$\infty$.  It is easy to check that
there exists a
(possibly empty) maximal open vertical strip in which the
integral~\eqref{mellin} is well
defined; that strip is called the \emph{fundamental strip}.

Suppose, for example, that~$f$ has the properties
\begin{equation}\label{tailbehavior}
f(x) = O(x^{b})\quad\hbox{as $x\to 0+$}\qquad \hbox{and}\qquad
f(x) = O(x^{a})\quad
\hbox{as $x\to \infty$}
\end{equation}
with $a < b$.
Let~$\alpha$ denote the infimum of such $a$'s, and $\beta$ the supremum of such
$b$'s.
Then $\langle - \beta, - \alpha \rangle$ is a substrip of the fundamental strip
[and equals the fundamental strip in typical cases, such as when
$f(x) = (1 + o(1)) c_0
x^{\beta}$ as $x \to 0+$ and $f(x) = (1 + o(1)) c_{\infty}
x^{\alpha}$ as $x \to \infty$ with
$c_0, c_{\infty} \neq 0$].
If, further, $f$ is continuously differentiable and monotone, then
one can check using
integration by parts that $\langle 1 - \beta, 1 - \alpha \rangle$ is
a substrip of the
fundamental strip for $f'$ and that
\begin{equation}
\label{mtdiff}
(f')^*(s) = - (s - 1) f^*(s - 1),\quad s \in \langle 1 - \beta, 1 -
\alpha \rangle.
\end{equation}

In this paper, we will make key use of the following Mellin transform pair:
\begin{equation}\label{logmt}
f(x) = \log(1 + x) \longleftrightarrow f^*(s) = \frac {\pi}{s \sin(\pi s)},
\end{equation}
with fundamental strip $\langle -1, 0 \rangle$.  For justification of
the double arrow, see
the start of the next subsection.

\subsection{Inversion and the mapping property}
\label{S:mtinv}

Suppose for simplicity that~$f$ is continuous and has a
nonempty fundamental strip, which we denote by
$\langle - \beta, - \alpha \rangle$.
Then for any~$c$ in the interval $( - \beta, - \alpha)$ we have
\begin{equation}\label{inversemellin}
f(x) = \frac 1{2 \pi i}\int_{c - i \infty}^{c + i \infty}\!f^*(s) x^{-s}\,ds.
\end{equation}
This establishes a correspondence between functions and their Mellin
transforms.

Mellin transforms are useful in asymptotic analysis because, under a certain
growth condition on the meromorphic continuation of the Mellin transform $f^*$,
there is a correspondence between
the asymptotic expansion of the function~$f$ near~$0$ or~$\infty$ and the
singularities of $f^*(s)$ for $s\in \mathbb C$.
(For our purposes, we need only asymptotic expansions near~$\infty$, so
those are all we will discuss here.)
Flajolet and Sedgewick~\cite{FS2003} call this correspondence the
\emph{mapping property}
of Mellin transforms.
To explain more fully, it is convenient to
introduce the following formal-sum shorthand notation.

Let $\varphi:\Omega \to {\mathbb C}$ be a meromorphic function in a
domain~$\Omega$ and let
${\cal P}\subset \Omega$ be the set of its poles.
For each $s_0 \in \Pc$,
let~$\Delta_{s_0}(s)$ be a truncation
of the Laurent series for $\varphi(s)$ at $s_0\in {\cal P}$; we
assume that $\Delta_{s_0}(s)$
contains at least all terms of order $(s - s_0)^{-k}$ with $k > 0$.
We then write
\begin{equation}
\label{asymp}
\varphi(s) \asymp \sum_{s_0 \in {\cal P}} \left[ \Delta_{s_0}(s)
\right]_{s = s_0} \quad (s \in \Omega),
\end{equation}
or, when no confusion results, simply
$$
\varphi(s) \asymp \sum_{s_0 \in {\cal P}} \Delta_{s_0}(s) \quad (s\in \Omega),
$$
and call \eqref{asymp} a \emph{singular expansion} of $\varphi$ in~$\Omega$.

The following two theorems are essential for our results.  The first of these
theorems gives (among other things) sufficient conditions for
meromorphic continuation
of the Mellin transform~$f^*$.
The second theorem shows that a singular expansion of $f^*$
characterizes the asymptotic behavior of~$f$.

\begin{theorem}[Direct Mapping Theorem]\label{t:dmt}
Let~$f$ have Mellin transform~$f^*$
at least in a nonempty strip
$\langle -\beta, -\alpha \rangle$.
Assume
that $f(x)$ admits as $x \to \infty$ a finite
asymptotic expansion
of the form
\begin{equation}
\label{expansion}
f(x) = \sum_{(\xi, k) \in A} c_{\xi, k}\,x^{\xi} (\log x)^k + O(x^{\gamma}),
\end{equation}
for a finite
set~$A$ of pairs $(\xi, k)$, where
the $\xi$'s satisfy
$\gamma < \xi \le \alpha$
and the $k$'s are nonnegative integers.
Then $f^*$ is continuable to a meromorphic function in the strip
$\langle -\beta, -\gamma \rangle$, where
it admits the singular expansion
$$
f^*(s) \asymp -\sum_{(\xi, k) \in A} c_{\xi, k} \frac {(-1)^k k!}{(s
+ \xi)^{k + 1}}\ \
(s \in \langle -\beta, -\gamma \rangle).~\nopf
$$
\end{theorem}

\begin{theorem}[Reverse Mapping Theorem]\label{t:rmt}
Let~$f$ be
continuous in $(0, \infty)$ with Mellin transform~$f^*$
existing at least in a nonempty strip
$\langle -\beta, -\alpha \rangle$.
\ignore{
\smallskip \\
(i)~Assume that~$f^*(s)$ admits a meromorphic continuation to the
strip $\langle -\gamma, -\alpha \rangle$
for some $\gamma > \beta$ with a finite number of poles there, and is analytic
on $\Re s = - \gamma$.  Assume also that there exists a real number
$\eta\in (\alpha, \beta)$
such that
\begin{equation}\label{growthcond}
f^*(s)=O(|s|^{-r}), \quad \hbox{with} \quad r>1,
\end{equation}
when $|s|\to \infty$ in $-\gamma \le \Re s \le -\eta$.  If~$f^*$
admits the singular
\marginal{Our~\eqref{sing} corrects a typo in~\cite{FS2003}:\ The
line before their (7.32)
should read $s \in \langle \gamma, \eta \rangle$.}
expansion
\begin{equation*}
\label{sing}
f^*(s) \asymp \sum_{(\xi, k) \in A} d_{\xi, k} \frac{1}{(s + \xi)^k}\ \
(s \in \langle -\gamma, -\eta \rangle)
\end{equation*}
with all $k$'s in the sum nonnegative integers,
then an asymptotic expansion
of $f(x)$ as $x \to 0+$ is
\begin{equation}
f(x) = \sum_{(\xi, k)\in A} d_{\xi,k} \left( \frac{(-1)^{k - 1}}{(k -
1)!} x^{\xi} (\log
x)^{k - 1} \right) + O(x^{\gamma}).
\end{equation}
}
Assume that $f^*(s)$ admits a meromorphic continuation
to
the strip
$\langle -\beta, -\gamma \rangle$
for some $\gamma < \alpha$
with a finite number of poles there,
and is analytic on $\Re s = - \gamma$.
Assume also that there exists a real number
$\eta \in (\alpha, \beta)$
such that
\begin{equation}\label{growthcond}
f^*(s) = O(|s|^{-r}), \quad\hbox{with}\quad r > 1,
\end{equation}
when $|s| \to \infty$ in $- \eta \le \Re s \le - \gamma$.
If~$f^*$ admits the singular expansion
$$
f^*(s) \asymp  - \sum_{(\xi, k)\in A} c_{\xi, k} \frac{(-1)^k k!}{(s
+ \xi)^{k + 1}}
\ \ (s \in \langle -\eta, -\gamma \rangle)
$$
with all $k$'s in the sum nonnegative integers,
then an asymptotic expansion
of $f(x)$ as $x \to \infty$ is
$$
f(x) = \sum_{(\xi, k)\in A} c_{\xi, k}\,x^{\xi} (\log x)^k +
O(x^{\gamma}).~\nopf
$$
\end{theorem}

\subsection{Termwise differentiation of an asymptotic expansion}
\label{S:termwise}

The proof of the following result is a simple exercise in the use of
the two mapping
theorems.  Note the condition $r > 2$.

\begin{corollary}
\label{C:termwise}
Let~$f$ be continuously differentiable and monotone, and let~$\alpha$
and~$\beta$ be defined
as in the paragraph containing~\eqref{tailbehavior}
and~\eqref{mtdiff}.  Assume that $\alpha <
\beta$, that $f(x)$ admits as $x \to \infty$ a finite asymptotic
expansion~\eqref{expansion}
as in the statement of the Direct Mapping Theorem, and that~$f^*$ has
a meromorphic
continuation that is analytic on the line $\Re s = - \gamma$ and
satisfies the growth
condition~\eqref{growthcond} with $r > 2$.  Then the
expansion~\eqref{expansion} can be
differentiated termwise:\ as $x \to \infty$,
$$
f'(x) = \sum_{(\xi, k) \in A} c_{\xi, k}\,x^{\xi - 1} (\log x)^{k -
1} (\xi \log x + k ) +
O(x^{\gamma - 1}).~\nopf
$$
\end{corollary}

\subsection{Harmonic sums and the separation property}
\label{S:harmonic}

The mapping property discussed in \refS{S:mtinv}
is particularly effective for asymptotic analysis of a harmonic sum
\begin{equation}\label{harmonic}
F(x) = \sum_k \lambda_k f(\mu_kx)
\end{equation}
with $\mu_k > 0$ for every~$k$.
The reason is the readily checked \emph{separation property}
that the Mellin transform of~\eqref{harmonic} has
the simple product form
\begin{equation}\label{mtharmonic}
F^*(s) = f^*(s) \sum_k \lambda_k\,\mu_k^{-s},
\end{equation}
valid for~$s$ in the intersection of the fundamental strip of~$f$ and
the domain of absolute
convergence of the \emph{generalized Dirichlet series} $\sum_k \gl_k
\mu_k^{-s}$.

\section{Asymptotic expansion for the log Laplace transform}
\label{S:hasy}
Recall
the expression~\eqref{hexpression} for the log Laplace transform $h_{\mm}
\equiv h$ of $V^2$.  To avoid nuisance factors it is convenient for us to
reparameterize~$h_{\mm}$ by defining, for $x\ge 0$,
\begin{equation}\label{lt}
L_{\mm}(x) := - 2 \log \EE\,\exp\left\{-\frac x2 V^2\right\}
            = \sum_{n_1 = 1}^{\infty}\!\!\cdots\!\!\sum_{n_d = 1}^{\infty}
                \log\left(1 + \frac{x}{\kappa_{n_1}(m_1) \cdots
                \kappa_{n_d}(m_d)}\right).
\end{equation}
Notice that $h_{\mm}(x) = \frac12 L_{\mm}(2x)$.
In this section we derive a complete asymptotic expansion for
$L_{\mm}$, for general~$\mm$.
First we treat the one-dimensional case.  Then we bring Mellin
asymptotic summation to bear.
That is, we factor the Mellin transform $L^*_{\mm}$ in terms of the transforms
$L^*_{m_j}$ using the separation property, analyze the singularities
of the factors using
(among other things) the Direct Mapping Theorem, and finally obtain
the asymptotic expansion
for $L_{\mm}$ (or for $h_{\mm}$) using the Reverse Mapping Theorem.
Using \refC{C:termwise}
we also derive corresponding expansions for the first two derivatives
of $h_{\mm}$, as needed
for our application in \refS{S:small} of Sytaja's Tauberian Theorem.

\subsection{The one-dimensional case}
\label{S:1d}

We begin by analyzing $m$-integrated Brownian motion, \ie, the $d = 1$ case
\begin{equation}\label{ellj}
L_{m}(x) = \sum_{n = 1}^{\infty} \log \left(1 + \frac x {\kappa_n(m)} \right).
\end{equation}
For integer~$\ell$, $0 \le \ell \le 2 m + 1$, set $\omega_{\ell} := \exp\{ i \pi
\frac{\ell}{m + 1}\} = \omega_1^{\ell}$, $v_{\ell} := \exp\{ i \pi \frac{2 \ell
+ 1}{2 m + 2}\} = v_0 \omega_{\ell}$, and $\beta_{\ell}(x) \equiv \beta_{\ell}
:= x^{1 / (2 m + 2)} i v_{\ell}$.  We employ the following result from Gao et
al.\ (cf.~Theorem~6 in~\cite{GaoHannigLeeTorcaso2003b}), which gives an exact
expression for $L_m(x)$:

\begin{lemma}\label{l:Lm}
For $x \geq 0$ we have
$L_m(x) = \log |\det N(x)| - (m + 1) \log(2 m + 2)$,
where
$$N(x) := \left(\begin{array}{llcl}1 & 1 & \cdots & 1\\
\omega_0 & \omega_1 & \cdots & \omega_{2m+1}\\
\vdots & \vdots & \cdots & \vdots \\
\omega_0^{m} & \omega_1^{m} & \cdots & \omega_{2m+1}^{m}\\
\omega_0^{m+1}e^{\beta_0} & \omega_1^{m+1}e^{\beta_1} & \cdots &
\omega_{2m+1}^{m+1}e^{\beta_{2m+1}}\\
\vdots & \vdots & \cdots & \vdots \\
\omega_0^{2m+1}e^{\beta_0} & \omega_1^{2m+1}e^{\beta_1} & \cdots &
\omega_{2m+1}^{2m+1}e^{\beta_{2m+1}}
\end{array}\right).~\nopf$$
\end{lemma}

The following important asymptotic consequence will be needed.

\begin{lemma}
\label{l:atinf}
As $x\to \infty$, the quantity $L_m(x)$ defined at~\eqref{ellj} has the
expansion
\begin{eqnarray*}
L_m(x)
   &=& \csc\left(\frac{\pi}{2 m + 2}\right) x^{1/(2 m + 2)} +  2 \log
|\det U| - (m + 1)
         \log(2 m + 2) \\
   & & \quad {} + O\left(\exp\left\{- \sin\left(\frac{\pi}{2 m +
2}\right) x^{1 / (2 m +
         2)}\right\}\right),
\end{eqnarray*}
where~$U_m \equiv U$ is the Vandermonde matrix
\begin{equation}
\label{Umatrix}
U := \left(\begin{array}{lllcl}1 & 1 & \cdots & 1\\
\omega_0 & \omega_1 & \cdots & \omega_{m}\\
\omega_0^2 & \omega_1^2 & \cdots & \omega_{m}^2\\
\vdots & \vdots & \cdots & \vdots \\
\omega_0^{m} & \omega_1^{m} & \cdots & \omega_{m}^{m}  \end{array}\right).
\end{equation}
\end{lemma}

\begin{proof}
We need to study the large-$x$ behavior of $N(x)$.
Multiply the last $m+1$ columns of $N(x)$ by $e^{\beta_0}$,
$e^{\beta_1}$, ..., $e^{\beta_m}$, respectively, and use
$\beta_j=-\beta_{m+1+j}$ to obtain the following matrix:
$$
\overline{N}(x) := \left(
\begin{array}{llllll}
1&\cdots&1&e^{\beta_{0}}&\cdots&e^{\beta_{m}}\\
\omega_0&\cdots&\omega_m&\omega_{m+1}e^{\beta_0}&\cdots
&\omega_{2m+1}e^{\beta_{m}}\\
\vdots&\cdots&\vdots&\vdots&\cdots&\vdots\\
\omega_0^m&\cdots&\omega_m^m&\omega_{m+1}^m
e^{\beta_0}&\cdots&\omega_{2m+1}^m e^{\beta_{m}}\\
\omega_0^{m+1}e^{\beta_0}&\cdots&\omega_m^{m+1}e^{\beta_m}&\omega_{m+1}^{m+1}
&\cdots&\omega_{2m+1}^{m+1}\\
\vdots&\cdots&\vdots&\vdots&\cdots&\vdots\\
\omega_0^{2m+1}e^{\beta_0}&\cdots&\omega_m^{2m+1}e^{\beta_m}&\omega_{m+1}^{2m+1}
&\cdots&\omega_{2m+1}^{2m+1}
\end{array}
\right).
$$
For $0\leq \ell \leq m$ we have
$$
|e^{\beta_{\ell}}| \leq |e^{\beta_0}| = \exp\{ - \sin\left( \pi / (2 m +
2) \right) x^{1 / (2 m + 2)}\}
$$
and therefore
$$\det \overline{N}(x)=\det \overline{N}(\infty) + O(|e^{\beta_0}|),
$$
where $\overline{N}(\infty)$ is the matrix obtained from
$\overline{N}(x)$ by replacing each entry with a factor $e^{\beta_{\ell}}$
by~$0$:
$$
\overline{N}(\infty) = \left(
\begin{array}{llllll}
1&\cdots&1&0&\cdots&0\\
\omega_0&\cdots&\omega_m&0&\cdots&0\\
\cdots&\cdots&\cdots&\cdots&\cdots&\cdots\\
\omega_0^m&\cdots&\omega_m^m&0&\cdots&0\\
0&\cdots&0&\omega_{m+1}^{m+1}&\cdots&\omega_{2m+1}^{m+1}\\
\cdots&\cdots&\cdots&\cdots&\cdots&\cdots\\
0&\cdots&0&\omega_{m+1}^{2m+1}&\cdots&\omega_{2m+1}^{2m+1}
\end{array}\right).
$$
A simple calculation reveals that $|\det \overline{N}(\infty)| =
|\det U|^2 \neq 0$,
where~$U$ is defined at~\eqref{Umatrix}.

Finally, using
$$
e^{\beta_0} e^{\beta_1} \cdots e^{\beta_m} =
     \exp\left\{ - \csc\left(\frac{\pi }{2m+ 2}\right) x^{1 / (2 m +
2)} \right\}
$$
we obtain
\begin{eqnarray*}
\lefteqn{L_m(x) + (m + 1) \log(2 m + 2)} \\
   &=& \log|e^{-\beta_0} e^{-\beta_1} \cdots e^{-\beta_m} \det
         \overline{N}(x)| \\
   &=& \log\left|\exp\left\{ \csc\left(\frac{\pi }{2m+ 2}\right) x^{1
/ (2 m + 2)}
         \right\} \det \overline{N}(x)\right| \\
   &=& \csc\left(\frac{\pi }{2m+ 2}\right) x^{1 / (2 m + 2)} +
         \log|\det \overline{N}(\infty)| + O(|e^{\beta_0}|) \\
   &=& \csc\left(\frac{\pi }{2m+ 2}\right) x^{1 / (2 m + 2)} + 2 \log|\det U| +
         O\hspace{-.04in}\left(\exp\left\{ - \sin\left(\frac{\pi }{2m+ 2}\right)
         x^{1 / (2 m + 2)} \right\}\right),
\end{eqnarray*}
as desired.
\end{proof}

\begin{remark}\label{R:fundstrip}
   From \refL{l:atinf} we see that $L_m(x) / x^{1/(2m+2)}$ has a finite
nonzero limit
as $x\to \infty$.  Also, it is easy to see that $L_m(x) / x$ has a
finite nonzero limit $x \to
0+$.  Therefore, the fundamental
strip for the Mellin transform $L_m^*$ of $L_m$ is
$\langle -1, -1/(2m+2)\rangle$.
\end{remark}

\begin{remark}\label{R:meromorphic}
As a consequence of
the Direct Mapping Theorem
we see that $L_m^*$ is
continuable
to a meromorphic function in the strip $\langle -1, \infty \rangle$, with
\begin{equation}\label{mtsingexp}
L_m^*(s) \asymp \left[ -\frac{\csc\left(\frac{\pi}{2 m + 2}\right)}
{s + \frac{1}{2 m + 2}} \right]_{s = - \frac{1}{2 m + 2}} + \left[
\frac {- 2 \log|\det U| +
(m + 1) \log(2 m + 2)}{s} \right]_{s = 0};
\end{equation}
in particular, in this strip $L_m^*(s)$ has only simple poles, at $s
= - 1 / (2 m + 2)$ and at
$s = 0$.
\end{remark}

\begin{remark}\label{R:vanderasy}
Although our treatment keeps~$m$ fixed, the reader may be curious as
to the large-$m$
behavior of the term $2 \log |\det U_m|$ appearing in \refL{l:atinf}.
Using the well-known
formula for the determinant of a Vandermonde matrix and approximating
the double sum that
results after taking logarithms by a double integral, we find
$$
2 \log |\det U_m| \sim m^2 \left[\frac{1}{2} \log 2 + \iint_{0 \leq x < y \leq
1}\log \sin\left(\frac{\pi}{2} (y - x)\right)\,dx\,dy\right] = -
\frac{7 \zeta(3)}{2 \pi^2}
m^2.
$$
We omit the details.
\end{remark}
\smallskip

Since $\displaystyle L_m(x)= \sum_{n = 1}^{\infty} \log
\left(1+\frac x{\kappa_n(m)}\right)$, it
follows from \eqref{logmt} and \eqref{mtharmonic} that
\begin{equation}\label{prelation}
L_m^*(s) = \frac {\pi}{s \sin(\pi s)} \sum_{n = 1}^{\infty} \kappa_n^s(m),
\end{equation}
or what is the same,
\begin{equation}\label{relation}
K_m(s) := \sum_{n = 1}^{\infty} \kappa_n^s(m) = \frac {s \sin(\pi
s)}{\pi} L_m^*(s).
\end{equation}
Since, as discussed in the next subsection, $\kappa_n(m)$ grows like
$n^{2 m + 2}$,
\eqref{prelation}--\eqref{relation} are valid (without continuation)
in the strip $\langle
-1, -1 / (2 m + 2) \rangle$.

\subsection{The generalized Dirichlet series $K_m(s)$}
\label{S:gds}

When we apply the Reverse Mapping Theorem to $L^*_{\mm}$ in
\refS{S:generalcase},
we will need to verify the growth condition~\eqref{growthcond}.  This
will rely on
corresponding estimates for the one-dimensional $L^*_m$, and so we
need to study the growth of
$L^*_m(s)$, or equivalently of the generalized Dirichlet series $K_m(s)$.

We can do this using the following very sharp asymptotic estimate for
$\kappa_n(m)$ (see also Gao et al.~\cite{GaoHannigLeeTorcaso2003a}).

\begin{lemma}[Gao et al.~\cite{GaoHannigTorcaso2003}, Theorem~2]
\label{L:kappaest}
For $m = 0$ we have $\kappa_n(0) = \left[ \left( n - \frac12 \right) \pi
\right]^2$ for every~$n$.  For each fixed $m \geq 1$, as $n \to \infty$ we have
$$
\kappa_n(m) = \left[ \left( n - \sfrac{1}{2} \right) \pi \right]^{2 m
+ 2} \left[ 1 + O\left( n^{-1} \exp\left\{ - \pi \sin\left( \frac{\pi}{m + 1}
\right) n \right\} \right) \right].~\nopf
$$
\end{lemma}
\noindent

 From \refL{L:kappaest} it is easily seen
that the domain of absolute convergence for $K_m$
is the strip
$\langle - \infty, - 1 / (2 m  + 2) \rangle$.  In this strip we can write
$$
K_m(s) = \Kh_m(s) + \sum_{n = 1}^{\infty} \left[ \kappa^s_n(m) - \kappah^s_n(m)
\right]
$$
where
$$
\Kh_m(s) := \sum_{n = 1}^{\infty} \kappah^s_n(m) = \left[ \left( \sfrac{\pi}{2}
\right)^{(2 m + 2) s} - \pi^{(2 m + 2) s} \right] \zeta( - (2 m + 2) s)
$$
with $\kappah_n(m) := \left[ \left( n - \sfrac{1}{2} \right) \pi \right]^{2 m +
2}$; here~$\zeta$ denotes Riemann's zeta function.  Using
Lemma \ref{L:kappaest} it is easy to check that the remainder series $\sum_{n =
1}^{\infty} \left[ \kappa^s_n(m) - \kappah^s_n(m) \right]$ converges absolutely
for all~$s \in {\mathbb C}$ and defines an entire function of~$s$
that is $O(|s|)$
in any strip $\langle - R, R \rangle$ with $0 < R < \infty$.  The generalized
Dirichlet series $\Kh_m(s)$ is meromorphically continuable for $s \in {\mathbb
C}$, with a single simple pole at $s = - 1 / (2 m + 2)$.  In any
strip $\langle -
R, R \rangle$, the continued $\Kh_m(s)$ grows at most polynomially
in~$|s|$ (see,
e.g.,~\cite{MR88c:11049}, Section~5.1).  Putting the pieces of the
above argument
together, we have established the following result.

\begin{lemma}
\label{L:gds}
For any $0 < R < \infty$, the meromorphic continuation of the
generalized Dirichlet series
$K_m(s)$ at~\eqref{relation} grows at most polynomially in $|s|$ as
$|s| \to \infty$ in the
strip $\langle - R, R \rangle$.~\nopf
\end{lemma}

\subsection{The general $d$-dimensional case:~separation}
\label{S:generalcase}

With
results regarding the one-dimensional processes in hand
we now turn our attention to arbitrary dimension~$d$ and vector~$\mm$. The
following simple lemma exploits the
harmonic-sum structure~\eqref{lt} of $L^*_{\mm}(s)$.
\begin{lemma}\label{l:induct}
For any dimension $d \ge 1$,
\begin{equation}\label{lmstarprod}
L_{\mm}^*(s) = \left(\frac{s \sin(\pi s)}{\pi}\right)^{d - 1}
\prod_{j = 1}^d L_{m_j}^*(s),
\end{equation}
valid (without continuation) in the strip $\langle - 1, - 1 / (2 m_1
+ 2) \rangle$.
\end{lemma}

\begin{proof}
When $d=1$ the result is trivial (recalling \refR{R:fundstrip}).
For $d \ge 2$, we observe
$$
L_{\mm}(x) = \sum_{n_d = 1}^{\infty} L_{(m_1, \dots, m_{d -
1})}\left(\frac{x}{\kappa_{n_d}(m_d)}\right),
$$
whence, by the separation property~\eqref{mtharmonic},
\eqref{relation}, and induction,
\begin{eqnarray*}
L_{\mm}^*(s)
   &=& L_{(m_1, \dots, m_{d - 1})}^*(s) \times \sum_{n_d = 1}^{\infty}
\kappa_{n_d}^s(m_d)
    =  \frac {s \sin(\pi s)}{\pi} L_{(m_1, \dots, m_{d - 1})}^*(s)
L^*_{m_d}(s) \\
   &=& \left(\frac{s \sin(\pi s)}{\pi}\right)^{d - 1} \prod_{j = 1}^d
L_{m_j}^*(s)
\end{eqnarray*}
for $s \in \langle - 1, - 1 / (2 m_1 + 2) \rangle$, as claimed.
\end{proof}

In preparation for our main theorem (\refT{t:main}) we next consider
the singular
expansion of the meromorphically continued $L^*_{\mm}(s)$.  In light of
\refL{l:induct} it is enough to understand the meromorphically continued
$L_{m_j}^*(s)$ for each $j$. Each
$L_{m_j}^*(s)$ has a simple pole at $s = 0$, but $s \sin(\pi s)$ has
a double zero at $s=0$;
thus, if $d \ge 2$,
then $L_{\mm}^*(s)$ has as singularities in $\langle -1, \infty \rangle$ only
poles at
$s = - 1 / (2 m_j + 2)$ for
$j = 1, \dots, d$.

Let
$$
r_m(x)
   := L_m(x) - \left[\csc\left(\frac{\pi}{2 m + 2}\right) x^{1/(2 m +
2)} +  2 \log |\det U| -
        (m + 1) \log(2 m + 2)\right]
$$
denote the exponentially small remainder term in the asymptotic
expansion for $L_m(x)$ in
\refL{l:atinf}.  From \refL{l:atinf} and (the proof of) the Direct
Mapping Theorem, for $s
\neq - 1 / (2 m_j + 2)$ we have
[compare~\eqref{mtsingexp}]
\begin{equation}\label{lstarasymp1}
L_{m_j}^*(s) = - \frac {\csc\left(\frac {\pi}{2 m_j + 2} \right)}
{s + \frac{1}{2 m_j + 2}} + A_{m_j}(s).
\end{equation}
Here
\begin{align}\label{Amjdef}
A_{m_j}(s)
&= \int_0^1\!L_{m_j}(x) x^{s - 1}\,dx + \frac{- 2 \log |\det U_{m_j}|
+ (m_j + 1)\log(2 m_j +
       2)}{s} \\
&  \qquad \qquad  \qquad {} + \int_1^{\infty}\!r_{m_j}(x) x^{s -
1}\,dx \nonumber
\end{align}
is analytic in the strip $\langle - 1, \infty \rangle$ except at $s = 0$.

\begin{remark}
\label{r:lmstarpoles}
 From the above discussion and \refL{l:induct} we immediately conclude
for $d \geq 2$ that
$L_{\mm}^*(s)$ is meromorphic in $\langle -1 , \infty \rangle$ with its only
poles located at  $s = - 1 / (2 m_j + 2)$ for $j=1,\dots ,d$.
Further, Lemmas~\ref{L:gds} and~\ref{l:induct}, together with the
exponential growth
of the sine function along vertical lines in~${\mathbb C}$, imply that
$L_{\mm}^*(s)$ satisfies the growth condition~\eqref{growthcond}
for \emph{every} $1 < r < \infty$.
\end{remark}

\subsection{Main theorem:\ Complete asymptotic expansion of $L_{\bf m}(x)$}
\label{S:asympLmx}

With the previous results we are now able to obtain a complete asymptotic
expansion for $L_{\bf m}(x)$. In order to state our main result, we suppose
\begin{equation*}
m_1 = \cdots = m_{t_1} < m_{t_1 + 1} = \cdots m_{t_1 + t_2} < \cdots < m_{t_1 +
\cdots + t_{g - 1} + 1} = \cdots = m_{t_1 + \cdots + t_g},
\end{equation*}
where the number of groups of ties is $g \ge 1$ and, for $1 \leq \nu \leq g$,
there is a tie of size $t_{\nu} \ge 1$ in group~$\nu$; thus $d = t_1 +
\cdots + t_g$. We will denote the common $m$-value in the $\nu$th group by
$\bar{m}_{\nu}$. Set
$$
{\xi}_{\nu} := \frac 1{2\bar{m}_{\nu}+2}.
$$
In order to use the Reverse Mapping Theorem, we need the singular expansion of
$L_{\mm}^*(s)$ near each pole $s = -{\xi}_{\nu}$ ($\nu = 1, \dots, g$).
From~\eqref{lmstarprod} and~\eqref{lstarasymp1}, near $s = -\xi_{\nu}$ we have
\begin{equation}
\label{lmstarfactored}
L_{\mm}^*(s) = \left(\frac{s \sin(\pi s)}{\pi}\right)^{d - 1}
\left(\prod_{n \ne \nu} \left[L_{\bar{m}_n}^*(s)\right]^{t_n}\right)
\left( \frac {-\csc\left(\pi \xi_{\nu} \right)}
{s + \xi_{\nu}} + A_{\bar{m}_{\nu}}(s)\right)^{t_{\nu}}.
\end{equation}
The product of the first two of the three factors on the right
in~\eqref{lmstarfactored} is analytic for~$s$ near $-\xi_{\nu}$
and therefore has a Taylor expansion that can be computed
up through the term of order $(s + \xi_{\nu})^{t_{\nu} - 1}$;
in regard to this computation, note that
repeated differentiation under the integrals in~\eqref{Amjdef} is
easily justified.
Similarly, we can expand the analytic function $A_{\bar{m}_{\nu}}(s)$
up through
order $(s + \xi_{\nu})^{t_{\nu} - 2}$
and use a mulitnomial
expansion to obtain a Laurent series for the third factor up through
order $(s + \xi_{\nu})^{-1}$.  Multiplying these expansions together
we can get a
Laurent series for~\eqref{lmstarfactored} up through
the term of
order $(s + \xi_{\nu})^{-1}$.  Applying the Reverse Mapping Theorem (the growth
condition is satisfied:\ see \refR{r:lmstarpoles}) we arrive at our
main result.

\begin{theorem}
\label{t:main}
For any $\mm$-integrated Brownian sheet in dimension $d \geq 2$,
$L_{\mm}(x)$ has
as $x \to \infty$ a complete asymptotic expansion of the following form:
\begin{equation}
\label{expansionform}
L_{\mm}(x) = \sum_{\nu=1}^g \sum_{k=0}^{t_{\nu}-1}
c_{\nu,k}\,x^{\xi_{\nu}}(\log
x)^k + O(x^{-R})
\end{equation}
for any $R > 0$.  The values of $c_{\nu,k}$ are computed as outlined in
the previous paragraph.
\nopf
\end{theorem}

\begin{remark}
\label{r:exposmall}
In fact, the error term in \refT{t:main} is exponentially small in a
positive power of~$x$.  We have written a complete proof, but will provide only
a sketch in the next paragraph (as the details are straightforward but
laborious), that the remainder term is bounded by $\exp\left\{ - c x^{1 / (d (2
m_d + 2))} \right\}$ for a certain constant~$c$ depending on ($d$ and) $\mm$.
We have not tried to optimize this bound; in particular, it may be that the
power $1 / (d (2 m_d + 2))$ can be improved to $1 / [2 (m_1 + \cdots + m_d) + 2
d]$.

The Reverse Mapping Theorem~\ref{t:rmt} is proved by invoking the
inverse Mellin
transform formula~\eqref{inversemellin} with $c = - \eta$ and then shifting the
line of integration rightward to $\Re s = - \gamma$ by means of the residue
theorem; the growth condition~\eqref{growthcond} is used to justify this shift
rigorously.  One then uses the growth condition again to bound the shifted
integral by $O(x^{\gamma})$.  When sharper growth estimates of $f^*(s)$ are
available, as they are in our case $f = L_{\mm}$ (cf.~\refR{r:lmstarpoles}),
one can let~$|\gamma|$ grow with~$x$ and obtain sharper remainder bounds.  Our
proof for $L_{\mm}$ uses well-known growth estimates for the Riemann
zeta function in vertical strips and takes $\gamma(x)  = - \ct\,x^{1 / (d (2 m_d
+ 2))} + O(1)$ for a suitable positive constant~$\ct$.
\end{remark}
\smallskip

Because $h(x) = \frac12 L_{\mm}(2x)$, it follows from \refT{t:main} that $h(x)$
also has an expansion of the form~\eqref{expansionform}, namely, for
any $R > 0$,
\begin{equation}
\label{hexpansion}
h(x) = \sum_{\nu=1}^g \sum_{k=0}^{t_{\nu}-1} c_{\nu,k}(0)\,x^{\xi_{\nu}}(\log
x)^k + O(x^{-R})
\end{equation}
where
\begin{equation*}
c_{\nu, k}(0) := 2^{\xi_{\nu} - 1} \sum_{\ell = k}^{t_{\nu} - 1} c_{\nu, \ell}
{\ell
\choose k} (\log 2)^{\ell - k}
\end{equation*}
(and error term as refined in \refR{r:exposmall}).  In order to apply Sytaja's
theorem in \refS{S:small}, we will need the following extension of \refT{t:main}
to the derivatives of~$h$.

\begin{lemma}
\label{l:hterm}
For $j = 0, 1, 2, \dots$, the function $h^{(j)}$ has as $x \to
\infty$ a complete
asymptotic expansion of the form
\begin{equation}
\label{hjexpansion}
h^{(j)}(x) = \sum_{\nu=1}^g \sum_{k=0}^{t_{\nu}-1} c_{\nu,k}(j)\,x^{\xi_{\nu}
- j}(\log x)^k + O(x^{-(R + j)})
\end{equation}
for any $0 < R < \infty$.  The expansions are obtained by successive termwise
differentiations, starting with~\eqref{hexpansion}.  In particular,
with $c_{\nu, k}(j) := 0$ whenever $k \geq t_{\nu}$,
\begin{align*}
c_{\nu, k}(1)
  &= \xi_{\nu} c_{\nu, k}(0) + (k + 1) c_{\nu, k + 1}(0), \\
c_{\nu, k}(2)
  &= (\xi_{\nu} - 1) c_{\nu, k}(1) + (k + 1) c_{\nu, k + 1}(1) \\
  &= \xi_{\nu} (\xi_{\nu} - 1) c_{\nu, k}(0) +  (2 \xi_{\nu} - 1) (k +
1) c_{\nu, k
       + 1}(0) + (k + 2) (k + 1) c_{\nu, k + 2}(0),
\end{align*}
and, for any $0 < R < \infty$,
\vspace{-.1in}
$$
h(x) - x h'(x) = \sum_{\nu=1}^g \sum_{k=0}^{t_{\nu}-1} [(1 - \xi_{\nu}) c_{\nu,
k}(0) - (k + 1) c_{\nu, k + 1}(0)] x^{\xi_{\nu}} (\log x)^k + O(x^{-R}).
$$
\end{lemma}

\begin{proof}
Denote the eigenvalues of~${\cal A}_{\mm}$ by $a_{\nn}$, as
at~\eqref{Vsum}.  Then
$h(x) = \sum_{\nn} \log(1 + a_{\nn} x)$ for all~$x \geq 0$, and the dominated
convergence theorem justifies the termwise differentiation giving,
for $j = 1, 2,
\dots$ and $x \geq 0$,
$$
h^{(j)}(x) = (-1)^{j - 1} (j - 1)! \sum_{\nn} a^j_{\nn} (1 + a_{\nn} x)^{- j}.
$$
In particular, each function~$h^{(j)}$ is monotone; and as $x \to 0+$
we have the
simple estimates
$$
h(x) = (1 + o(1))\,x \sum_{\nn} a_{\nn}, \qquad h^{(j)}(x) = (1 +
o(1)) (-1)^{j -
1} (j - 1)! \sum_{\nn} a^j_{\nn}.
$$
It is then easy to use induction on~$j$, together with \refC{C:termwise} and
\refR{r:lmstarpoles}, to complete the proof.
\end{proof}

\begin{remark}
\label{r:jexposmall}
For each $j = 0, 1, 2, \dots$, the error term in~\eqref{hjexpansion} is
exponentially small in a positive power of~$x$; compare \refR{r:exposmall}.
\end{remark}

\subsection{Expansion of $L_{\bf m}(x)$:\ examples}
\label{S:asympLmxexamples}

We next give three simple examples of the computations entering into
the expansion
in \refT{t:main}.

\subsubsection{Lead-order asymptotics in the general $d$-dimensional
case ($d \geq 2$)}
\label{S:leadordergeneric}

Suppose $d \geq 2$ and
\begin{equation*}
m := m_1=\cdots =m_t<m_{t+1}\le \cdots \le m_d.
\end{equation*}
According to \refT{t:main}, the lead-order asymptotics for
$L_{\mm}(x)$ are given
by
\begin{equation*}
L_{\mm}(x) =  (1 + o(1))\,c_{1, t - 1}\,x^{1 / (2 m + 2)}(\log x)^{t-1},
\end{equation*}
since (as will become evident below) $c_{1, t - 1} > 0$.
This lead-order term corresponds via the Reverse Mapping Theorem to
the term of order $(s + \xi_1)^{- t}$ in the expansion of
$L^*_{\mm}(s)$ near its
pole at $s = - \xi_1$.
From~\eqref{lmstarfactored} we see immediately that, as $s \to -\xi_1$,
\begin{align*}
L_{\mm}^*(s)
 &= (1 + o(1))
      \left( \frac{\xi_1 \sin(\pi \xi_1)}{\pi} \right)^{d - 1}
      \left( \prod_{j = t  +1}^d L_{m_j}^*(-\xi_1) \right)
      \left( \frac {-\csc\left(\pi \xi_1 \right)} {s + \xi_1} \right)^t \\
  &= - (1 + o(1))
      \frac{\left[ \sin\!\left( \frac {\pi}{2 m + 2} \right) \right]^{d - 1 -
      t}}{[(2 m + 2) \pi]^{d - 1}(t - 1)!}
      \left( \prod_{j = t + 1}^d L_{m_j}^*\!\left(\!-\frac 1{2m+2}\right)
      \right) \frac {(-1)^{t - 1}(t - 1)!} {(s + \xi_1)^t}
\end{align*}
from which it follows
that
\begin{equation}
\label{leadc}
c_{1, t - 1} = \frac{\left[ \sin\!\left( \frac {\pi}{2 m + 2} \right) \right]^{d
- 1 - t}}{[(2 m + 2)\pi]^{d - 1}(t - 1)!} \left( \prod_{j = t + 1}^d
L_{m_j}^*\!\left(\!-\frac{1}{2 m + 2}\right) \right).
\end{equation}

Thus, in general, computation of the lead-order term for $L_{\mm}(x)$ requires
(numerical) evaluation of the fundamental-strip Mellin transform values
\begin{equation}
\label{Mellinintegration}
L^*_{m_j}\!\left( - \frac{1}{2 m + 2} \right) = \int^{\infty}_0\!L_{m_j}(x)\,x^{
- (2 m + 3) / (2 m + 2)}\,dx > 0.
\end{equation}
However, in the particular equal-$m$'s case where $t = d$ we obtain the simpler
result
$$
L_{\mm}(x) = (1 + o(1)) \frac{\csc(\frac{\pi}{2 m + 2})}{[(2 m + 2) \pi]^{d -
1}(d - 1)!} x^{1 / (2 m + 2)}(\log x)^{d - 1}.
$$

\begin{remark}
\label{R:russiandifference}
 From~\eqref{relation} we see that an alternative to numerical evaluation of
the integral~\eqref{Mellinintegration} is numerical evaluation of the
generalized Dirichlet series $K_{m_j}\left(- \frac{1}{2 m + 2}\right)$.
Evaluation of such sums arises in the approach
of~\cite{KarolNazarovNikitin2003}, which is more computationally intensive
since it involves numerical root-finding.
\end{remark}

\subsubsection{Complete asymptotic expansion: distinct $m_j$'s}
\label{S:distinctexpexample}

Given \refT{t:main}, it is almost trivial to obtain a full asymptotic
expansion for $L_{\mm}(x)$ when the vector~$\mm$ has distinct components.
That is, suppose that $g = d \geq 2$, so that $t_{\nu} = 1$ for $\nu
= 1, \dots,
d$.  Then \refT{t:main} implies
\begin{equation}
\label{Ldgeq2}
L_{\mm}(x) = \sum_{\nu = 1}^d C_{\nu}\,x^{1 / (2 m_{\nu} + 2)} + O(x^{-R})
\end{equation}
for any $0 < R < \infty$, where
\begin{equation}
\label{Cnu}
C_{\nu} := \frac{\left[\sin\!\left(\frac{\pi}{2 m_{\nu} +
2}\right)\right]^{d - 2}
\prod_{j \neq {\nu}} L^*_{m_j}\!\left(-\frac{1}{2 m_{\nu} +
2}\right)}{[( 2 m_{\nu}
+ 2) \pi]^{d - 1}}.
\end{equation}

\subsubsection{Full asymptotic expansion: $d=2$}
\label{S:fullexpexample}

\refT{t:main} provides a complete asymptotic expansion of $L_{\mm}(x)$ for
any vector~$\mm$.  However, calculations quickly become quite
cumbersome to perform
by hand; symbolic-manipulation software such as {\tt Mathematica} or
{\tt Maple}
can then be of great help.  Complementing the above result for
distinct $m$'s, here
is another case where we can spell out the expansion without too much notation.

Suppose that $m := m_1 = m_2$, i.e.,\ that $g = 1$ and $d = t_1 = 2$.  We then
find
$$
L_{\mm}(x) = c_{1, 1}\,x^{1 / (2 m + 2)} \log x + c_{1, 0}\,x^{1 / (2 m + 2)} +
O(x^{-R})
$$
for any $0 < R < \infty$, where
\begin{align}
\label{constequalm1}
c_{1, 1} &= \frac{\csc\!\left(\frac{\pi}{2 m + 2}\right)}{(2 m + 2) \pi}, \\
\label{constequalm2}
c_{1, 0} &= \frac{A_m\!\left(- \frac{1}{2 m + 2}\right)}{(m + 1) \pi} +
\frac{\csc\!\left(\frac{\pi}{2 m + 2}\right)}{\pi} +
\frac{\cos\!\left(\frac{\pi}{2 m + 2}\right) \csc^2\!\left(\frac{\pi}{2 m +
2}\right)}{2 m + 2},
\end{align}
and $A_m(\cdot)$ is computed via~\eqref{Amjdef}.  When $m = 0$ in this
example, i.e.,\ when~$V$ is the $L^2$-norm of \emph{two-dimensional
Brownian sheet},
we have from \refL{l:Lm} that
$$
L_0(x) \equiv \log \cosh(x^{1 / 2})
$$
and
$$
A_0(- 1 / 2) = \int_0^1\!L_0(x)\,x^{- 3 / 2}\,dx - 2 \log 2 +
\int_1^{\infty}\!\left[L_0(x) - x^{1 / 2} + \log 2\right] x^{- 3 /
2}\,dx \approx -0.3624.
$$
Thus,
$$
L_{(0, 0)}(x) = \frac{1}{2 \pi} x^{1 / 2} \log x + c_{1, 0}\,x^{1 /
2} + O(x^{-R})
$$
for any $0 < R < \infty$, and here
$$
c_{1, 0} = \frac{1 + A_0( - 1 / 2)}{\pi} \approx 0.2029.
$$

\section{Reversion of $h'(x) = \eps$}
\label{S:reversion}

In principle we are now in position to use Sytaja's theorem in conjunction with
the asymptotic expansions for $h(x) - x h'(x)$ and $h''(x)$ in \refL{l:hterm} to
state a strong small-deviations result for integrated Brownian sheets.  The
problem is that, for more explicit results, we still must solve $h'(x) =
\eps$ to get $x^* \equiv x^*(\eps)$ in order to obtain the needed expansions
in~$\eps$, up to an additive term $o(1)$ for $h(x^*) - \eps\,x^*$ and up to a
factor $1 + o(1)$ for $h''(x^*)$.

In this \refS{S:reversion} we discuss the reversion of the asymptotic expansion
for $h'(x)$ given in \refL{l:hterm}.  We begin in \refS{S:asysuff} by showing
that the exponentially small remainders in all our asymptotic expansions may
safely be ignored when carrying out reversion and applying Sytaja's theorem.
In \refS{S:leadorderreversion} we determine the lead-order asymptotics
for~$x^*$.  Clearly, the lead-order asymptotics agree for~$x^*$ and the
(unique real) solution $x_0 \equiv x_0(\eps)$ to the equation
\begin{equation}
\label{x0eq}
\eps = \sum_{k = 0}^{t_1 - 1} c_{1, k}(1)\,x^{\xi_1 - 1}(\log x)^k
\end{equation}
obtained from the equation $\eps = h'(x)$ by including only those terms with
$\nu = 1$ (corresponding to the largest power of~$x$) in the
expression~\eqref{hjexpansion} (with $j = 1$) for $h'(x)$.  Next, in
\refS{S:reversiongeneral} we obtain a complete asymptotic expansion for
$x^*$ in terms of elementary functions and $x_0$.  Finally, in \refS{S:x0}
we discuss exact computation of~$x_0$.

\subsection{Truncating the asymptotic expansion for~$h'$ suffices}
\label{S:asysuff}

Our next result implies that, as far as obtaining lead-order asymptotics for
the small-ball probability, we may act as if the asymptotic expansions in
\refL{l:hterm} were exact expressions.

\begin{lemma}
\label{l:asysuff}
Let~$\hh$ denote the truncated asymptotic expansion for~$h$:
\begin{equation*}
\hh(x) := \sum_{\nu = 1}^g \sum_{k = 0}^{t_{\nu} - 1} c_{\nu,
k}(0)\,x^{\xi_{\nu}} (\log x)^k,
\end{equation*}
and let~$\xh$ satisfy $\hh(\xh)=\varepsilon$.
As $\eps \to 0+$, the following errors all tend to zero faster than any
power of~$\eps$:
$$
\xh - x^*,\quad \hh^{(j)}(\xh) - h^{(j)}(x^*)\mbox{ (for $j = 0, 1,
\dots$)},\quad [\hh(\xh) - \xh \hh'(\xh)] - [h(x^*) - \eps\,x^*].
$$
\end{lemma}
Observe that $\hh'(x)$ decreases monotonically to~$0$ for
sufficiently large~$x$;
thus for sufficiently small~$\eps$ there exists a unique solution $\xh \equiv
\xh(\eps)$ to $\hh'(\xh) = \eps$.
We remark in passing that $\hh^{(j)}$ can also be obtained by truncating the
asymptotic expansion for $h^{(j)}$ given in \refL{l:hterm}.

\begin{proof}
 From definitions and \refL{l:hterm} we have
$$
\hh'(\xh) = \eps = h'(x^*) = \hh'(x^*) + O((x^*)^{-R})\mbox{\ \ for any $0 <
R < \infty$},
$$
which by the simple \refL{l:leadonceandforall} below can be written
$$
\hh'(\xh) - \hh'(x^*) = O(\eps^R)\mbox{\ \ for any $0 < R < \infty$},
$$
as $\eps \to 0+$.  It follows using the mean value theorem that
$$
\xh - x^* = O(\eps^R)\mbox{\ \ for any $0 < R < \infty$}.
$$
The other assertions follow readily.
\end{proof}

\begin{remark}
The errors in \refL{l:asysuff} are each exponentially small in a positive
power of $1 / \eps$.  This assertion follows easily from the case $j = 1$ of
\refR{r:jexposmall}.
\end{remark}

\subsection{Lead-order reversion in the general $d$-dimensional
case ($d \geq 2$)}
\label{S:leadorderreversion}

If all we want is a standard ``weak'' small-deviations result, we need only
obtain the lead-order asymptotics for $x^*(\eps)$
[equivalently, for $x_0(\eps)$].  This is easy.  Recall from \refL{l:hterm} that
$h'(x) = (1 + o(1))\,c_{1, t - 1}(1)\, x^{\xi - 1}\,(\log x)^{t - 1}$, where $t
:= t_1$ is the number of $m_j$'s equal to the smallest value $m := m_1$ and
$\xi := \xi_1 = 1 / (2 m + 2)$.  The coefficient~$c_{1, t - 1}(1)$ is given by
$c_{1, t - 1}(1) = \xi 2^{\xi - 1} c_{1, t - 1}$, with $c_{1, t - 1}$ given
by~\eqref{leadc}.  The proof of the following lemma is left to the reader.

\begin{lemma}
\label{l:leadonceandforall}
If $h'(x^*) = \varepsilon$ and~$x_0$ is the solution to~\eqref{x0eq}, then
as $\varepsilon \to 0+$ we have $x^* = (1 + o(1))\,x_0 = (1 + o(1))\,\xt_0$,
where
\begin{equation}
\label{leadx}
\xt_0 \equiv \xt_0(\eps) := \left[\frac{c_{1, t - 1}(1)}{(1 - \xi)^{t - 1}}
\cdot \frac{1}{\eps} \left(\log \frac{1}{\varepsilon}\right)^{t - 1}
\right]^{\frac{1}{1 - \xi}}.~\qed
\end{equation}
\end{lemma}

We will use \refL{l:leadonceandforall} in
\refS{S:small} to obtain a
rather simple ``weak'' small-deviations result for arbitrary~$\mm$; see
\refT{t:weaksmallball}.
\refL{l:leadonceandforall} also provides us with all the asymptotic information
about $h''(x^*)$ we need in applying Sytaja's theorem:

\begin{corollary}
\label{c:hdp}
If $h'(x^*) = \eps$, then, with~$\xt_0$ defined at~\eqref{leadx}, as $\eps
\to 0+$ we have
$$
- (x^*)^2 h''(x^*) = (1 + o(1))\,(1 - \xi)\,\xt_0\,\eps
$$
\end{corollary}

\begin{proof}
This is routine.  Recall \refL{l:hterm} and in particular the coefficients
$$
c_{1, t - 1}(2) = \xi (\xi - 1) 2^{\xi - 1} c_{1, t - 1} = - (1 - \xi) c_{1, t
- 1}(1).
$$
Then
\begin{align*}
- (x^*)^2 h''(x^*)
 &= - (1 + o(1))\,\xt_0^2\,h''(\xt_0) = - (1 + o(1))\,c_{1, t -
      1}(2)\,\xt_0^{\xi}\,(\log \xt_0)^{t - 1} \\
 &= - (1 + o(1))\,c_{1, t - 1}(2)\,\xt_0^{\xi}\,\left( \frac{1}{1 - \xi} \log
      \frac{1}{\eps} \right)^{t - 1},
\end{align*}
and it is easy to check from definitions that
$$
- c_{1, t - 1}(2)\,\xt_0^{\xi}\,\left( \frac{1}{1 - \xi} \log \frac{1}{\eps}
\right)^{t - 1} =  (1 - \xi) \xt_0 \eps.~\qed
$$
\noqed
\end{proof}

\subsection{Reversion: an asymptotic expansion for~$x^*$}
\label{S:reversiongeneral}

We next derive a complete asymptotic expansion for~$x^*$ in terms of
elementary functions and the solution~$x_0$ to~\eqref{x0eq}.  We will discuss
exact computation of~$x_0$ in \refS{S:x0}.

Dropping the error term from the expansion~\eqref{hjexpansion} (with $j = 1$)
for $h'(x)$ (as justified by \refL{l:asysuff}), let us write the equation $\eps
= h'(x)$ in the form
\vspace{-.1in}
\begin{equation}
\label{epsf}
\eps = \sum_{\nu = 1}^g f_{\nu}(x),
\end{equation}
where
\begin{equation}
\label{fnudef}
f_{\nu}(x) := x^{- \eta_{\nu}} \sum_{k = 0}^{t_{\nu} - 1} a_{\nu, k} (\log x)^k,
\quad \nu = 1, \dots, g;
\end{equation}
here, for abbreviation, we have set $a_{\nu, k} := c_{\nu, k}(1)$ and
$\eta_{\nu} := 1 - \xi_{\nu} = (2 \bar{m}_{\nu} + 1) / (2 \bar{m}_{\nu} + 2)$,
so that
$$
1 / 2 \leq \eta_1 < \eta_2 < \cdots \eta_g < 1.
$$
(The notation $\bar{m}_{\nu}$ is as in \refS{S:asympLmx}.)  If $g = 1$, then we
have simply $x^* = x_0$; so we assume now that $g \geq 2$.

The main result of this subsection is the following complete asymptotic
expansion for $x^* / x_0$ in terms of inductively defined quantities $y_j$
($j = 0, 1, 2, \dots$).  Let $y_0 := 1$.  Suppose that $j \geq 1$ and that we
have defined $y_0, \dots, y_{j - 1}$.  Then set
$$
y^+_{j - 1} := y_0 + y_1 + \cdots + y_{j - 1} \quad \mbox{and} \quad
x_{j - 1}   := x_0\,y^+_{j - 1}
$$
and define
\begin{equation}
\label{yjdef}
y_j := \frac{y^+_{j-1} \left[ \sum_{\nu = 1}^g f_{\nu}(x_{j - 1}) - \eps
\right]}{\eta_1 f_1(x_{j - 1}) - x_{j - 1}^{- \eta_1} \sum_{k = 0}^{t_1 - 1} k
a_{1, k} (\log x_{j - 1})^{k - 1}}.
\end{equation}

\begin{proposition}
\label{p:reversiongeneral}
For each $j = 0, 1, \dots$ we have
\begin{equation}
\label{xjexpansion}
\frac{x^*}{x_0} = 1 + y_1 + y_2 + \cdots + (1+o(1))y_j
\end{equation}
and
\begin{equation}
\label{yjest}
y_j = O\!\left( \left(\eps^{\frac{\eta_2}{\eta_1} - 1} \left( \log
\frac{1}{\eps} \right)^{\Delta} \right)^j \right),
\end{equation}
where $\Delta := \max\{ t_{\nu} - 1 - (t_1 - 1) \frac{\eta_{\nu}}{\eta_1}: 2
\leq \nu \leq g \}$.
\end{proposition}

Note that, with $x_{- 1} := 0$, \eqref{xjexpansion} can be written equivalently
as
\begin{equation}
\label{xjalt}
x^* = x_{j - 1} + (1 + o(1)) x_0 y_j \quad (j = 0, 1, \dots).
\end{equation}
Before we prove the proposition, we illustrate its application in two special
cases.

\begin{example}[$t_1 = 1$]
\label{E:t1=1}
Suppose $t_1 = 1$, whence $f_1(x)$ reduces to $a_{1, 0} x^{- \eta_1}$ and it is
elementary that $x_0 = (a_{1, 0} / \eps)^{1 / \eta_1}$.  In this case,
\eqref{yjdef} simplifies to
\begin{equation}
\label{simpleryj}
y_j = \frac{y^+_{j-1}}{\eta_1 a_{1, 0}} x_{j - 1}^{\eta_1} \left[ \sum_{\nu =
1}^g f_{\nu}(x_{j - 1}) - \eps \right].
\end{equation}
If we assume further, as in \refS{S:distinctexpexample}, that all the $m_j$'s
are distinct (\ie,\ $g = d$), then $f_{\nu}(x)$ simplifies to $a_{\nu,
0}\,x^{- \eta_{\nu}}$ and~\eqref{simpleryj} can be written in the form
\begin{equation}
\label{yjdistinct}
y_j = \frac{(y^+_{j-1})^{1 + \eta_1}}{\eta_1} \left[ \sum_{\nu = 1}^d a_{\nu,
0}\,a_{1, 0}^{-
\eta_{\nu} / \eta_1} \eps^{(\eta_{\nu} / \eta_1) - 1} (y^+_{j - 1})^{-
\eta_{\nu}} - 1 \right].
\end{equation}
 From this it is easy to prove by induction that each $y_j$ has an asymptotic
expansion in increasing powers of~$\eps$ wherein each power is a nonnegative
integer combination of the numbers $(\eta_{\nu} / \eta_1) - 1$ with $2 \leq \nu
\leq d$.  The same is therefore true of $x^* / x_0$.
\end{example}

\begin{example}[$m_1 = 0$]
\label{E:m1=0}
Suppose $m_1 = \bar{m}_1 = 0$, whence $\eta_1 = 1 / 2$; since $\bar{m}_2 \geq
1$, we also have $\eta_2 \geq 3 / 4$.  In this case it is not hard to see
that the finite expansion
\begin{align*}
\frac{x^*}{x_0}
 &= 1 + y_1 + y_2 + O(y_3) = 1 + y_1 + y_2 + O\!\left( \eps^{6 \eta_2 - 3}
      \left( \log \frac{1}{\eps} \right)^{3 \Delta} \right) \\
 &= 1 + y_1 + y_2 + O\!\left( \eps^{3 / 2} \left( \log \frac{1}{\eps} \right)^{3
        \Delta} \right)
\end{align*}is sufficient
to obtain the lead-order asymptotics for the small-ball probability.
To be concrete, suppose further
that there are no ties, so that we are in the context of~\eqref{yjdistinct}
and $x_0 = (a_{1, 0} / \eps)^2$.  Then
\begin{align}
y_1 &= 2 \sum_{\nu = 2}^d a_{\nu, 0}\,a_{1, 0}^{- 2 \eta_{\nu}} \eps^{2
         \eta_{\nu} - 1}, \nonumber \\
y_2 &= - 4 \left[ \sum_{\nu = 2}^d a_{\nu, 0}\,a_{1, 0}^{- 2 \eta_{\nu}} \eps^{2
         \eta_{\nu} - 1} \right] \left[ \sum_{\nu = 2}^d \left( \eta_{\nu} -
         \frac34 \right) a_{\nu, 0}\,a_{1, 0}^{- 2 \eta_{\nu}} \eps^{2
         \eta_{\nu} - 1} \right] + O(\eps^{6 \eta_2 - 3}) \nonumber \\
\label{y2}
    &= - 4 \left( \eta_2 - \frac34 \right) a_{2, 0}^2\,a_{1, 0}^{- 4
         \eta_2} \eps^{4 \eta_2 - 2} + O(\eps^{2 (\eta_2 + \eta_3) - 2}) +
         O(\eps^{6 \eta_2 - 3}),
\end{align}
where the first remainder term in~\eqref{y2} is used if $d \geq 3$, and the
second if $d = 2$.  So, whether $\eta_2 > 3 / 4$ or $\eta_2 = 3 / 4$, and
whether $d \geq 3$ or $d = 2$,
$$
\frac{x^*}{x_0} = 1 + y_1 + o(\eps) = 1 + 2 \sum_{\nu = 2}^d a_{\nu, 0}\,a_{1,
0}^{- 2 \eta_{\nu}} \eps^{2 \eta_{\nu} - 1} + o(\eps),
$$
which we will see in \refE{E:small;distinct;m1=0} is sufficient
to obtain the lead-order asymptotics for the small-ball probability.
\end{example}
\medskip

\begin{proof}[Proof of \refP{p:reversiongeneral}]
We prove~\eqref{xjalt} and~\eqref{yjest} together by induction on~$j$.  The
base case $j = 0$ of the induction is  simple.  For $j \geq 1$, our induction
hypothesis is that the $(j - 1)$st instances of~\eqref{xjalt} and~\eqref{yjest}
hold; in particular, we can write
$$
x^* = x_{j - 2} + x_0 y_{j - 1} + x_0 y = x_{j - 1} + x_0 y = x_{j - 1} \left( 1
+ \frac{y}{y^+_{j - 1}} \right)
$$
for some $y = o(y_{j - 1})$.  For $\nu = 1, \dots, g$ we then have,
from~\eqref{fnudef},
\begin{align*}
f_{\nu}(x^*)
 &= f_{\nu}\!\left( x_{j - 1} \left( 1 + \frac{y}{y^+_{j - 1}} \right) \right)
      \\
 &= x_{j - 1}^{- \eta_{\nu}} \left( 1 + \frac{y}{y^+_{j - 1}} \right)^{-
      \eta_{\nu}}\,\sum_{k = 0}^{t_{\nu} - 1} a_{\nu, k} \left[ \log x_{j - 1} +
      \log \left( 1 + \frac{y}{y^+_{j - 1}} \right) \right]^k.
\end{align*}
For $\nu \geq 2$ we conclude
\begin{align}
f_{\nu}(x^*) - f_{\nu}(x_{j - 1})
 &= O(y\,f_{\nu}(x_{j - 1})) = O(y\,x_0^{- \eta_{\nu}} (\log x_0)^{t_{\nu} -
      1}) \nonumber \\
\label{2est}
 &= O\!\left( y\,\eps^{\eta_{\nu} / \eta_1} \left( \log \frac{1}{\eps}
      \right)^{t_{\nu} - 1 - (t_1 - 1) (\eta_{\nu} / \eta_1)} \right)
  = O\!\left( y\,\eps^{\eta_2 / \eta_1} \left( \log \frac{1}{\eps}
      \right)^{\Delta} \right),
\end{align}
where the third equality follows from \refL{l:leadonceandforall}.  For $\nu =
1$ we find
\begin{equation}
\label{1est}
f_1(x^*) - f_1(x_{j - 1})
 = - (1 + o(1)) \frac{y}{y^+_{j - 1}} \left[ \eta_1 f_1(x_{j - 1}) -
x_{j - 1}^{- \eta_1} \sum_{k = 0}^{t_1 - 1} k a_{1, k} (\log x_{j - 1})^{k-1}
\right],
\end{equation}
which we note is of the same order as $y \eps$.  Therefore,
$$
\eps = \sum_{\nu = 1}^g f_{\nu}(x^*) = \sum_{\nu = 1}^g f_{\nu}(x_{j - 1}) - (1
+ o(1)) \frac{y}{y^+_{j - 1}} \left[ \eta_1 f_1(x_{j - 1}) - x_{j - 1}^{-
\eta_1} \sum_{k = 0}^{t_1 - 1} k a_{1, k} (\log x_{j - 1})^{k-1} \right],
$$
and so [recalling~\eqref{yjdef}] $y = (1 + o(1)) y_j$.  This establishes the
$j$th instance of~\eqref{xjalt}.

To establish the $j$th instance of~\eqref{yjest}, we note again that the
denominator of~\eqref{yjdef} is precisely of order~$\eps$.  We must also
estimate the difference appearing in the numerator.  Suppose that $j
\geq 2$; the estimation for $j = 1$ is similar and left to the reader.  By
slight modifications of the arguments for~\eqref{2est}--\eqref{1est}, for $\nu
\geq 2$ we have
\begin{equation}
\label{diffnu}
f_{\nu}(x_{j - 1}) - f_{\nu}(x_{j - 2})
  = O\!\left( y_{j - 1}\,\eps^{\eta_2 / \eta_1} \left( \log \frac{1}{\eps}
      \right)^{\Delta} \right),
\end{equation}
and, also utilizing~\eqref{yjdef}, for $\nu = 1$ we have
\begin{align}
f_1(x_{j - 1}) - f_1(x_{j - 2})
 &= - \frac{y_{j - 1}}{y^+_{j - 2}} \left[ \eta_1 f_1(x_{j - 2}) - x_{j - 2}^{-
        \eta_1} \sum_{k = 0}^{t_1 - 1} k a_{1, k} (\log x_{j - 2})^{k-1} \right] +
        O\!\left( y_{j - 1}^2 \eps \right) \nonumber \\
\label{diff1}
 &= - \left[ \sum_{\nu = 1}^g f_{\nu}(x_{j - 2}) - \eps \right] +
        O\!\left( y_{j - 1}^2 \eps \right).
\end{align}
Summing the~$g$ equations \eqref{diffnu}--\eqref{diff1} and rearranging,
$$
\sum_{\nu = 1}^g f_{\nu}(x_{j - 1}) - \eps = O\!\left( y_{j -
1}\,\eps^{\eta_2 / \eta_1} \left( \log \frac{1}{\eps} \right)^{\Delta} \right) +
O\!\left( y_{j - 1}^2 \eps \right).
$$
By the induction hypothesis and our assumption that $j \geq 2$, the first
of the two $O(\cdot)$ terms predominates and the estimate
$$
\sum_{\nu = 1}^g f_{\nu}(x_{j - 1}) - \eps =
O\!\left( \eps \left(\eps^{\frac{\eta_2}{\eta_1} - 1} \left( \log
\frac{1}{\eps} \right)^{\Delta} \right)^j \right)
$$
is established, completing the induction.
\end{proof}

\subsection{Computation of~$x_0$}
\label{S:x0}

There remains the task of solving the following equation to obtain~$x_0$:
\begin{equation}
\label{ef1}
\eps = f_1(x) = \sum_{k = 0}^{t - 1} a_k\,x^{- \eta} (\log x)^k,
\end{equation}
with $t := t_1$, $a_k := a_{1, k} = c_{1, k}(1)$, and $\eta := \eta_1 = 1 -
\xi_1 = (2 m_1 + 1) / (2 m_1 + 2)$.  As has already been noted in \refE{E:t1=1},
this is trivial if $t = 1$:\ then $x_0 = (a_0 / \eps)^{1 / \eta}$.  For
general~$t$, of course, one can resort (for given $\eps > 0$) to Newton's
method or other root-finding methods.  What we will derive in this subsection
is a series representation of the solution, from which will follow an
asymptotic expansion for~$x_0$ sufficient to yield a complete asymptotic
expansion for the logarithm of the small-ball probability, but \emph{not} for
that probability itself (for which the exact value of~$x_0$ should be used).

It is instructive to begin by considering the case $t = 2$.  Setting $w := -
\eta [\log x +(a_0 / a_1)]$, equation~\eqref{ef1} can be rewritten as
$$
w e^w = - \frac{\eta}{a_1} \exp\{- \eta a_0 / a_1\} \eps.
$$
This equation for~$w$ cannot be solved in terms of elementary functions, but it
does have the solution
$$
w = W\!\left( - \frac{\eta}{a_1} \exp\{- \eta a_0 / a_1\} \eps \right),
$$
in terms of the Lambert $W$-function.  (See~\cite{CorlessEtal1996}; more
specifically, $W$ here is the branch $W_{-1}$ in the notation there.  We could
alternatively work in terms of the closely related glog function defined
in~\cite{Kalman1996}.)  According to the paragraph following equation~(4.20)
in~\cite{CorlessEtal1996} we have, in terms of $z := (\eta / a_1) \exp\{- \eta
a_0 / a_1\} \eps$, the expansion
\begin{equation}
\label{weq1}
w = W(- z) = \log z - \log \log \frac{1}{z} + \sum_{r = 0}^{\infty} \sum_{s =
1}^{\infty} d_{r s} \left( \log \log \frac{1}{z} \right)^s (\log z)^{- (r +
s)}
\end{equation}
with $d_{r s}$ given in terms of Stirling numbers of the first kind (see
equation~1.2.9-(26) in~\cite{knuth97}) by $d_{r s} = \frac{1}{s!}
(-1)^r \stirlingone{r + s}{r + 1}$.  As remarked in~\cite{CorlessEtal1996},
the series in~\eqref{weq1} is absolutely and uniformly convergent for
small~$\eps$ and also serves as an
asymptotic expansion for $w - \log z + \log \log \frac{1}{z}$ as $\eps \to 0+$.
It is not hard to show that this result
can be rearranged to one of the following form, wherein the series has the same
convergence and asymptotic expansion properties as in~\eqref{weq1}:
\begin{equation}
\label{weq2}
x_0 = \xt_0 \times \left[ 1 + \sum_{r = 0}^{\infty} \sum_{s = 0}^{\infty} \dt_{r
s} \left( \log \log \frac{1}{\eps} \right)^s \left( \log \frac{1}{\eps}
\right)^{- (r + s)} \right]\mbox{\ \ with $\dt_{0 0} = 0$}.
\end{equation}
Recall that, in general,
\begin{equation}
\label{xt0}
\xt_0 := \left[ \frac{a_{t - 1}}{\eta^{t - 1}} \cdot \frac{1}{\eps} \cdot
\left(\log \frac{1}{\eps}\right)^{t - 1} \right]^{1 / \eta}
\end{equation}
is the lead-order approximant to~$x_0$ found in \refL{l:leadonceandforall}, and
for $t = 2$ this reduces to $\xt_0 = \left(\frac{a_1}{\eta} \cdot
\frac{1}{\eps} \cdot \log \frac{1}{\eps}\right)^{1 / \eta}$.  In particular, the
first correction term in~\eqref{weq2} then has coefficient $\dt_{0 1} = d_{0 1}
/ \eta = 1 / \eta$.

The expansion~\eqref{weq2} can be extended to general values of~$t$:

\begin{lemma}
\label{l:x0expansion}
For general~$t$, the solution~$x_0$ to~\eqref{ef1} has an expansion of the form
\eqref{weq2}--\eqref{xt0}.  The series in~\eqref{weq2} is absolutely and
uniformly convergent for small $\eps > 0$, and when rearranged to
\begin{equation}
\label{series}
\sum_{\ell = 1}^{\infty} \sum_{s = 0}^{\ell} \dt_{\ell - s, s} \left( \log \log
\frac{1}{\eps} \right)^s \left( \log \frac{1}{\eps} \right)^{- \ell}
\end{equation}
provides a complete asymptotic expansion for $(x_0 / \xt_0) - 1$.  The lead term
in~\eqref{series} has coefficient $\dt_{0 1} = (t - 1)^2 / \eta$.
\end{lemma}

\begin{proof}
We give only a sketch of the proof.  We can write the result of
\refL{l:leadonceandforall} in the form
\begin{equation}
\label{wgeneq}
w := - \eta \log x_0 = - \sigma^{-1} - (t - 1) \sigma^{-1} \tau - \lambda
- v,
\end{equation}
with $\lambda := \log(a_{t - 1} / \eta^{t - 1})$, where $v = o(1)$
and
$$
\sigma := \frac{1}{\log(1 / \eps)} \quad \mbox{and} \quad
\tau   := \frac{\log \log(1 / \eps)}{\log(1 / \eps)}.
$$
Substituting~\eqref{wgeneq} into the defining equation~\eqref{ef1} for~$x_0$,
we find
$$
\eps
  = e^w \sum_{k = 0}^{t - 1} a_k\,( - w / \eta)^k
  = e^w \sum_{k = 0}^{t - 1} a_k\,\eta^{ - k} \left[ \sigma^{-1} + (t - 1)
\sigma^{-1} \tau + \lambda + v \right]^k,
$$
or equivalently that~$v$ is a root of the function
\begin{align*}
F(\zeta)
 &:= e^{\zeta} - 1 - \left\{ \left[ 1 + (t - 1) \tau + \lambda \sigma +
\sigma \zeta \right]^{t - 1} - 1 \right\} \\
 &\qquad \quad {} - e^{- \lambda} \sigma^{t - 1} \sum_{k
= 0}^{t - 2} a_k\,\eta^{ - k} \sigma^{-k} \left[ 1 + (t - 1) \tau + \lambda
\sigma + \sigma \zeta \right]^k.
\end{align*}

The proof now proceeds as on pp.\ 347--349 in~\cite{CorlessEtal1996} to obtain
a series representation for~$v$, and then we exponentiate to
obtain~\eqref{weq2}.  We omit the details and will be content here to verify the
asserted value of $\dt_{0 1}$.  Indeed, our proof sketch shows that
$\eta\,\dt_{0 1}$ equals
$$
\frac{1}{2 \pi i} \int_{|\zeta| = \pi}\!G(\zeta)\,d\zeta,
$$
where the integration is taken counterclockwise and $G(\zeta)$ is the
coefficient of $\sigma^0 \tau^1$ in $\zeta F'(\zeta) / F(\zeta)$, namely,
$$
G(\zeta) = (t - 1)^2 \zeta e^{\zeta} \left( e^{\zeta} - 1 \right)^{-2}.
$$
Evaluating the integral completes the proof.
\end{proof}

\begin{remark}
In \refS{S:small} we will see that, when viewed as an asymptotic expansion for
$x_0$,  \eqref{weq2} is not accurate enough to use in obtaining a strong
small-deviations result but is good enough (see \refT{t:weakcorr}) to yield an
expansion for the logarithm of the small-ball probability.
\end{remark}

\section{Small deviation estimates}
\label{S:small}

We now have at hand all of the ingredients we need to apply Sytaja's theorem
and obtain $p \equiv p(\eps) := \Pr(V^2 \leq \eps)$ up to a factor $1 + o(1)$,
\ie,\ its log up to an additive $o(1)$.  Until \refS{S:d=1} we will assume $d
\geq 2$.  (The results for $d = 1$ are simpler to derive but require some
modification.)

To see how we are now positioned to obtain $- \log p$ up to additive $o(1)$,
note by Sytaja's theorem and \refC{c:hdp} that
\begin{equation}
\label{logp}
- \log p = h(x^*) - \eps x^* + \mbox{$\frac12$} \log \left[ 2 \pi (1 - \xi_1)
\eps
\xt_0 \right] + o(1)
\end{equation}
where $\xi_1 = 1 / (2 m_1 + 2)$ and $\xt_0$ is given explicitly
by~\eqref{leadx}.  Here~$x^*$ solves $h'(x) = \eps$, and
\refP{p:reversiongeneral} gives an asymptotic expansion allowing its
computation to an arbitrarily large power of~$\eps$ once the solution~$x_0$
to~\eqref{x0eq} is obtained (as discussed in \refS{S:x0}).  The expansion
for~$x^*$ can then be substituted into the expansion~\eqref{hexpansion} for
$h(\cdot)$.

In this section we present for the reader's convenience some explicit
small-deviations estimates in a few cases where the final result of the above
program is reasonably clean.  In \refS{S:weak} we determine lead-order
asymptotics, and more, for $- \log p$.  In \refS{S:strong} we determine
lead-order asymptotics for~$p$ itself.  For completeness, in \refS{S:d=1} we
handle the one-dimensional case of $m$-integrated Brownian motion.

\subsection{Logarithmic small-deviations estimates ($d \geq 2$)}
\label{S:weak}

Our first two results give lead-order asymptotics for $ - \log p$.

\begin{theorem}\label{t:weaksmallball}
Let $V$ be given as at~\eqref{Vdef} for an arbitrary $(m_1, \dots,
m_d)$-integrated Brownian sheet with $d \geq 2$, and assume $m := m_1 = \cdots =
m_t < m_{t+1} \le \cdots \le m_d$.  Then as $\varepsilon \to 0+$ we have
\vspace{-.2in}
\begin{equation*}
-\log \Pr\left(V^2 \le \varepsilon\right) = (1 + o(1)) \frac
{2 m + 1}{2} \left[\frac{C (2 m + 2)^{t - 2}}{(2 m + 1)^{t - 1}}
\right]^{\!\!\frac{2 m + 2}{2 m + 1}}
\left(\!\frac 1{\eps}\!\right)^{\!\!\frac{1}{2 m + 1}}
\!\left(\!\log\frac 1{\eps}\!\right)^{\!\!\!\!\frac{(t - 1)(2 m + 2)}{2 m
+ 1}},
\end{equation*}
where
\vspace{-.1in}
\begin{equation}
\label{cdef}
C(d, t, m) \equiv C := \frac{\left[\sin\left(\frac{\pi}{2 m +
2}\right)\right]^{d - 1 - t}
\prod_{j = t + 1}^d L_{m_j}^*\left(-\frac{1}{2 m + 2}\right)}{(t -
1)!\,[(2m + 2) \pi]^{d - 1}}.
\end{equation}
\end{theorem}

\begin{proof}
Using Sytaja's Tauberian Theorem and \refL{l:hterm} we find
$$
- \log \Pr(V^2 \leq \eps) = (1 + o(1))\,[h(x^*) - \eps\,x^*] = (1 + o(1)) 2^{\xi
- 1} (1 - \xi) c_{1, t - 1}\,(x^*)^{\xi} (\log x^*)^{t - 1},
$$
where $\xi = 1 / (2 m + 2)$ and $c_{1, t - 1}$ is given by~\eqref{leadc} and so
equals~$C$.  Now use \refL{l:leadonceandforall} and rearrange to obtain the
desired result.
\end{proof}

\begin{corollary}\label{C:d=tleadorder}
If $m_1 = \cdots = m_d = m$,
then as $\varepsilon\to 0+$ we have
\begin{eqnarray*}
\lefteqn{-\log \Pr\left(V^2 \le \varepsilon\right)} \\
  &=& (1 + o(1)) \frac {2m+1}2
\left[\frac {(2m+2)^{-1}\csc\!\left(\frac {\pi}{2m+2}\right)}
{(d-1)!\,[(2m+1)\pi]^{d-1}}\right]^{\frac {2m+2}{2m+1}}
\!\!\!\left(\!\frac 1{\varepsilon}\!\right)^{\!\!\frac 1{2m+1}}
\!\left(\!\log \frac 1{\varepsilon}\!\right)^{\!\!\!\!\frac
{(d-1)(2m+2)}{2m+1}}.~\qed
\end{eqnarray*}
\end{corollary}
\refC{C:d=tleadorder} was previously established by Li~\cite{Li1992a}.  More
recently, \refT{t:weaksmallball} was obtained independently of us, and
essentially simultaneously, by Karol' et.~al~\cite{KarolNazarovNikitin2003}
with a different form for the leading constant.  Indeed, the trace sums that
appear in Corollary~5.2 of~\cite{KarolNazarovNikitin2003} are simply generalized
Dirichlet series values:\ recall our \refR{R:russiandifference}.

Next we show that it is possible to obtain arbitrarily high-order corrections to
the estimate for $ - \log p$ in \refT{t:weaksmallball}.  [Note, however, that
the asymptotic scale clearly is not fine enough to obtain $- \log p$ up to
additive $o(1)$.]

\begin{theorem}
\label{t:weakcorr}
The $o(1)$ expression in \refT{t:weaksmallball} has a complete asymptotic
expansion of the form
$$
\sum_{r = 1}^{\infty} \sum_{s = 0}^r D_{r s} \left( \log \log
\frac{1}{\eps} \right)^s \left( \log \frac{1}{\eps} \right)^{- r}.
$$
The lead-order coefficient is $D_{1 1} = (t - 1)^2 (2 m + 2) / (2 m + 1)$.
\end{theorem}

\begin{proof}
We give a sketch.  Using \refP{p:reversiongeneral} with $j = 1$, one finds that
substitution of~$x_0$ for~$x^*$ in~\eqref{logp} produces an error which is
negligible relative to the scale involved in the statement of the theorem.
Substitution of the expansion in \refL{l:x0expansion} for~$x_0$ leads to the
desired result.  The details are straightforward but quite tedious, and are
omitted.
\end{proof}

\refT{t:weakcorr} immediately gives us the following two examples.

\begin{example}
When $\mm = {\bf 0}$, \refT{t:weakcorr} sharpens Csaki's classical
result~\cite{Csaki1982} for $\mbox{$d$-dimensional}$ Brownian sheet ($d \geq
2$) to
\begin{align*}
-\log \Pr(V^2 \le \eps)
 &= \frac{1}{8} \left[ (d - 1)!\,\pi^{d - 1} \right]^{-2}
      \frac{1}{\eps} \left(\log \frac{1}{\eps} \right)^{2 (d - 1)} \\
 &{} \qquad \times \left[ 1 + 2 (d - 1)^2\,\frac{\log \log \frac{1}{\eps}}{\log
       \frac{1}{\eps}} + O\!\left( \frac{1}{\log \frac{1}{\eps}} \right)
       \right].
\end{align*}
\end{example}

\begin{example}
\label{e:equalm2} For the $(m,m)$-integrated Brownian sheet, $m\ge
0$, we have
\begin{eqnarray*}
-\log \Pr (V^2\le \eps )
&=&
  \alpha \left(\frac{1}{\eps}\right)^{\frac{1}{2 m + 1}} \left(\log
    \frac{1}{\eps}\right)^{\frac{2 m + 2}{2 m + 1}} + \alpha \cdot \frac{2 m +
    2}{2 m + 1} \left(\frac{1}{\eps} \log \frac{1}{\eps}\right)^{\frac{1}{2 m +
    1}} \log \log \frac{1}{\eps} \\
& &
  \qquad {} + O\!\left(\left(\frac{1}{\eps} \log
    \frac{1}{\eps}\right)^{\frac{1}{2 m + 1}}\right),
\end{eqnarray*}
where
$$
\alpha(m) \equiv \alpha := \frac{1}{2} (2 m + 1)^{- \frac{1}{2 m + 1}} \left[
\frac{\csc\!\left( \frac{\pi}{2 m + 2} \right)}{(2 m + 2) \pi} \right]^{\frac{2
m + 2}{2 m + 1}}.
$$
\end{example}

\subsection{Strong small-deviation estimates ($d \geq 2$)}
\label{S:strong}

We focus on the case of distinct $m_j$'s.  In this case ($t = 1$) the statement
of \refT{t:weaksmallball} simplifies to $ - \log \Pr\left(V^2 \le \eps\right) =
(1 + o(1)) E(\eps)$, where~$C$ is given by~\eqref{cdef} and $E(\eps)$ is the
expression
\begin{equation}
\label{Edef}
E(\eps) = \frac{2 m + 1}{2} \left( \frac{C}{2 m + 2} \right)^{\!\!\frac{2 m +
2}{2 m + 1}} \left( \frac{1}{\eps} \right)^{\!\!\frac{1}{2 m + 1}}.
\end{equation}

\begin{theorem}
\label{t:distinctstrong}
Let $V$ be given as at~\eqref{Vdef} for an $(m_1, \dots, m_d)$-integrated
Brownian sheet with $d \geq 2$, and suppose that $m := m_1 < \cdots < m_d$. Then
as $\varepsilon \to 0+$ we have
$$
\Pr\left(V^2 \le \varepsilon\right) = (1 + o(1)) \left[ \frac{\pi}{m + 1}
E(\eps) \right]^{- 1 / 2} \exp\{ - E(\eps) \times [1 + \Sigma(\eps)] \}
$$
for some finite linear combination~$\Sigma(\eps)$ of powers of~$\eps$, wherein
each power is a nonzero nonnegative integer combination of the numbers
$$
\frac{(2 m_1 + 2) (2 m_{\nu} + 1)}{(2 m_1 + 1) (2 m_{\nu} + 2)} - 1, \qquad \nu
= 2, \dots, d.
$$
\end{theorem}

\begin{proof}
This follows routinely from~\eqref{logp} and \refE{E:t1=1}, recalling the
notation $\eta_{\nu} = (2 m_{\nu} + 1) / (2 m_{\nu} + 2)$ for $\nu = 1, \dots,
d$.  We omit the details.
\end{proof}

\begin{example}
\label{E:d=2,distinctm}
Even when $d = 2$, some vectors~$\mm$ require many terms in $\Sigma(\eps)$ in
order to obtain the small-ball probability up to a factor $1 + o(1)$.  We will
be content here to illustrate \refT{t:distinctstrong} by using it to state a
first-order correction to \refT{t:weaksmallball} when $d = 2$ and $m_1 < m_2$.
The present example thus complements \refE{e:equalm2}.

For the $(m_1,m_2)$-integrated Brownian sheet with $m_1 < m_2$, we have
$$-\log \Pr (V^2\le \eps ) = \alpha_1
\left(\frac 1{\eps}\right)^{\!\!\frac {1}{2m_1+1}}
+ \alpha_2 \left(\frac 1{\eps}\right)^{\!\!\frac {1}{2m_1+1}\cdot
\frac {2m_1+2}{2m_2+2}} +\,
o\left(\!\!\left(\frac 1{\eps}\right)^{\!\!\frac {1}{2m_1+1}\cdot
\frac {2m_1+2}{2m_2+2}}\!\right),
$$
where
$\alpha_1 = (2m_1+1)[c_{1,0}(1)]^{\frac {2m_1+2}{2m_1+1}}$, $\alpha_2 =
\left(2m_2 + 2\right)[c_{1,0}(1)]^{\frac {1}{2m_1+1}\cdot
\frac {2m_1+2}{2m_2+2}}c_{2,0}(1)$, and $c_{1,0}(1)$ and $c_{2,0}(1)$ are given as
in \refL{l:hterm}.
\end{example}

\begin{example}
\label{E:small;distinct;m1=0}
Suppose, as in \refE{E:m1=0}, that $m = m_1 = 0$ and that the $m_j$'s are
distinct, in which case~\eqref{cdef} reduces to
$$
C = (2 \pi)^{- (d - 1)} \prod_{j = 2}^d L_{m_j}^*(- 1 / 2)
$$
and~\eqref{Edef} simplifies further to
$$
E(\eps) = \frac{C^2}{8} \times \frac{1}{\eps}.
$$
Then, using the calculations in~\refE{E:m1=0} it is easy to check that the only
powers of~$\eps$ appearing in
$\Sigma(\eps)$ of \refT{t:distinctstrong} are powers $2 \eta_{\nu} - 1$, for
$\nu = 2, \dots, d$, and $4 \eta_2 - 2$; we recall $\eta_{\nu} := (2 m_{\nu} +
1) / (2 m_{\nu} + 2)$.  If, in particular, $d = 2$, then only two terms are
needed in $\Sigma(\eps)$:\ see the next example.
\end{example}

\begin{example}
For $(0, m_2)$-integrated Brownian sheet with $m_2 > 0$, we have
$$
\Pr(V^2 \le \eps) = (1 + o(1)) \frac{\eps^{1 / 2}}{c_{1,0}(1) \sqrt{\pi}}
\exp\left\{ d_1 \cdot \frac{1}{\eps} + d_2\!\left( \frac{1}{\eps}
\right)^{\!\!\frac{1}{m_2+1}} + d_3\!\left( \frac{1}{\eps} \right)^{\!\!- \frac
{m_2 - 1}{m_2 + 1}} \right\}
$$
where
$d_1 = [c_{1,0}(1)]^2$, $d_2 = (2m_2+2)[c_{1,0}(1)]^{\frac{1}{m_2 + 1}}c_{2,0}$,
and
$d_3 = [c_{1,0}(1)]^{- \frac{2 m_2}{m_2 + 1}} [c_{2,0}(1)]^2$, and~$c_{1,0}(1)$
and~$c_{2,0}(1)$ are given as in \refL{l:hterm}.
\end{example}

\subsection{Strong small deviation estimates when $d = 1$}
\label{S:d=1}
The difference between the cases $d \geq 2$ and $d = 1$ can be seen by
comparing~\eqref{Ldgeq2} and \refL{l:atinf} and noting the constant term $2 \log
|\det U_m| - (m + 1) \log(2 m + 2)$ in \refL{l:atinf}.  Once this difference is
noted, it is simple to obtain results for $d = 1$ like those in the preceding
two subsections.  In particular, the logarithmic small-deviations result in
\refC{C:d=tleadorder} remains true verbatim when $d = 1$.  Moreover, we have
the following strong small-deviations analogue, studied
previously in~\cite{GaoHannigTorcaso2003}, \cite{GaoHannigLeeTorcaso2003a},
\cite{GaoHannigLeeTorcaso2003b}, \cite{Nazarov2003},
and~\cite{NazarovNikitin2003}, of \refT{t:distinctstrong}:

\begin{theorem}
\label{t:d=1strong}
Let $V^2 \equiv V_m^2 := \int_0^1\!X^2_m(t)\,dt$, where~$X_m$ is
$m$-times integrated Brownian motion defined in \refR{R:mintBM}.
Then as $\varepsilon \to 0+$ we have
$$
\Pr\left(V^2 \le \varepsilon\right) = (1 + o(1)) \left[ | \det U_m | \cdot (2 m
+ 2)^{(m + 1) / 2} \right]^{-1} \left[ \frac{\pi}{m + 1} E(\eps) \right]^{- 1 /
2} \exp\{ - E(\eps) \}.
$$
Here, with $C = \csc(\pi / (2 m + 2))$, the expression $E(\eps)$ is given
by~\eqref{Edef}, and~$U_m$ is the Vandermonde matrix defined at~\eqref{Umatrix}.
\end{theorem}

\bibliographystyle{plain}
\begin{small}

\end{small}
\end{document}